\documentclass{amsart}
\usepackage{times,amssymb,amscd,graphics,axodraw,verbatim,youngtab}
\usepackage[all]{xy}

\numberwithin{equation}{section}

\newtheorem{prop}{Proposition}
\newtheorem{theorem}[prop]{Theorem}

\newtheorem{lemma}[prop]{Lemma}

\theoremstyle{definition}

\newtheorem{example}[prop]{Example}

\numberwithin{prop}{section}

\newcommand{\ld}{\dot{\ell}}
\newcommand{\ldt}{\tilde{\dot{\ell}}}
\newcommand{\ltd}{\dot{\tilde{\ell}}}
\newcommand{\s}{s}
\newcommand{\ldb}{\dot{s}}
\newcommand{\ltb}{\tilde{s}}
\newcommand{\ldtb}{\tilde{\dot{s}}}
\newcommand{\ltdb}{\dot{\tilde{s}}}

\newcommand{\nut}{\tilde{\nu}^\bullet}
\newcommand{\nud}{\dot{\nu}^\bullet}
\newcommand{\nutd}{\dot{\tilde{\nu}}^\bullet}
\newcommand{\nudt}{\tilde{\dot{\nu}}^\bullet}
\newcommand{\Jtil}{\tilde{J}^\bullet}
\newcommand{\Jd}{\dot{J}^\bullet}
\newcommand{\Jtd}{\dot{\tilde{J}}^\bullet}
\newcommand{\Jdt}{\tilde{\dot{J}}^\bullet}

\newcommand{\A}{\mathcal{A}}
\newcommand{\bh}{\hat{b}}
\newcommand{\bj}{\mathrm{bj}}
\newcommand{\bs}{\mathrm{bs}}
\newcommand{\bt}{\tilde{b}}
\newcommand{\btt}{\tilde{\bt}}
\newcommand{\Bbh}{\hat{\bar{B}}}
\newcommand{\Bh}{\hat{B}}
\newcommand{\Bt}{\widetilde{B}}
\newcommand{\cc}{cc}
\newcommand{\comp}{\theta}
\newcommand{\Conf}{\mathrm{C}}
\newcommand{\D}{D}
\newcommand{\Dh}{\hat{D}}
\newcommand{\dist}{\mathrm{dist}}
\newcommand{\ds}{\delta_s}
\newcommand{\dspin}{-_s}
\newcommand{\dt}{\tilde{\delta}}
\newcommand{\Dt}{\tilde{D}}
\newcommand{\dual}{*}
\newcommand{\emb}{\mathrm{emb}}
\newcommand{\es}{\varnothing}
\newcommand{\F}{F}
\newcommand{\Fh}{\hat{F}}
\newcommand{\Ft}{\tilde{\F}}
\newcommand{\gggg}{\mathfrak{g}}
\newcommand{\id}{\mathrm{id}}
\newcommand{\inner}[2]{\langle #1\,,\,#2\rangle}
\newcommand{\Jht}{\hat{J}^\bullet}
\newcommand{\Jt}{J^\bullet}
\newcommand{\Jtt}{\tilde{J}^\bullet}
\newcommand{\Jttt}{\tilde{\tilde{J}}^\bullet}
\newcommand{\La}{\Lambda}
\newcommand{\la}{\lambda}
\newcommand{\lah}{\hat{\la}}

\newcommand{\lb}{\overline{\ell}}
\newcommand{\Lbh}{\hat{\bar{L}}}
\newcommand{\Lh}{\hat{L}}
\newcommand{\lh}{\mathrm{lh}}
\newcommand{\lhs}{\lh_s}
\newcommand{\lbt}{\tilde{\lb}}

\newcommand{\lt}{\tilde{\ell}}

\newcommand{\mt}{\tilde{m}}
\newcommand{\nb}{\overline{n}}
\newcommand{\nht}{\hat{\nu}^\bullet}
\newcommand{\nt}{\nu^\bullet}
\newcommand{\ntt}{\tilde{\nu}^\bullet}
\newcommand{\nttt}{\tilde{\tilde{\nu}}^\bullet}
\newcommand{\ol}{\overline}
\newcommand{\om}{\omega}
\newcommand{\omb}{\overline{\omega}}
\newcommand{\Path}{\mathcal{P}}
\newcommand{\Phit}{\tilde{\Phi}}
\newcommand{\qbin}[2]{\genfrac{[}{]}{0pt}{}{#1}{#2}}
\newcommand{\RC}{\mathrm{RC}}
\newcommand{\rh}{\mathrm{rh}}
\newcommand{\sbar}{\overline{s}}
\newcommand{\sbt}{\tilde{\sbar}}

\newcommand{\st}{\tilde{s}}

\newcommand{\sq}{\Box}
\newcommand{\tj}{\mathrm{tj}}
\newcommand{\ts}{\mathrm{ts}}
\newcommand{\vdom}{\scalebox{0.5}{\yng(1,1)}}
\newcommand{\vL}{\boldsymbol{L}}

\newcommand{\wt}{\mathrm{wt}}
\newcommand{\Z}{\mathbb{Z}}

\begin{document}

\pagestyle{plain}

\title{A bijection between type $D_n^{(1)}$ crystals and rigged configurations}

\author[A.~Schilling]{Anne Schilling}
\address{Department of Mathematics, University of California, One Shields
Avenue, Davis, CA 95616-8633, U.S.A.}
\email{anne@math.ucdavis.edu}
\urladdr{http://www.math.ucdavis.edu/\~{}anne}
\thanks{\textit{Date:} June 2004}
\thanks{Partially supported by the NSF grant DMS-0200774.}

\begin{abstract}
Hatayama et al. conjectured fermionic formulas associated with tensor products 
of $U'_q(\mathfrak{g})$-crystals $B^{r,s}$. The crystals $B^{r,s}$ correspond
to the Kirillov--Reshetikhin modules which are certain finite dimensional
$U'_q(\mathfrak{g})$-modules.
In this paper we present a combinatorial description of the affine
crystals $B^{r,1}$ of type $D_n^{(1)}$. A statistic preserving bijection between
crystal paths for these crystals and rigged configurations is given,
thereby proving the fermionic formula in this case. This bijection reflects 
two different methods to solve lattice models in statistical mechanics: the
corner-transfer-matrix method and the Bethe Ansatz.
\end{abstract}

\maketitle

\section{Introduction}
The quantized universal enveloping algebra $U_q(\gggg)$ associated
with a symmetrizable Kac--Moody Lie algebra $\gggg$ was introduced
independently by Drinfeld \cite{D:1985} and Jimbo \cite{J:1985} in their
study of two dimensional solvable lattice models in statistical
mechanics. The parameter $q$ corresponds to the temperature of the
underlying model. Kashiwara \cite{K:1990} showed that at zero
temperature or $q=0$ the representations of $U_q(\gggg)$ have
bases, which he coined crystal bases, with a beautiful
combinatorial structure and favorable properties such as
uniqueness and stability under tensor products.

The irreducible finite-dimensional $U'_q(\gggg)$-modules for an affine 
Kac--Moody algebra $\gggg$ were classified by Chari and 
Pressley~\cite{CP:1995,CP:1998} in terms of Drinfeld polynomials. The 
Kirillov--Reshetikhin modules $W^{r,s}$, labeled by a Dynkin node $r$ 
(different from the $0$-node) and a positive integer $s$, form a special 
class of these finite-dimensional modules. They naturally correspond
to the weight $s\La_r$, where $\La_r$ is the $r$-th fundamental weight
of $\gggg$. Recently, Hatayama et al.~\cite{HKOTY:1999,HKOTT:2002}
conjectured that the Kirillov--Reshetikhin modules $W^{r,s}$ have
a crystal basis denoted by $B^{r,s}$.
The existence of such crystals allows the definition of one dimensional 
configuration sums, which play an important role in the study of phase 
transitions of two dimensional exactly solvable lattice models. 
For $\gggg$ of type $A^{(1)}_n$, the existence of the crystal $B^{r,s}$ was
settled in \cite{KKMMNN:1992}, and the one dimensional configuration sums 
contain the Kostka polynomials, which arise in the theory of symmetric 
functions, combinatorics, the study of subgroups of finite abelian groups, and 
Kazhdan--Lusztig theory. In certain limits they are branching functions of
integrable highest weight modules.

In~\cite{HKOTY:1999,HKOTT:2002} fermionic formulas for the one dimensional
configuration sums were conjectured. Fermionic formulas originate in
the Bethe Ansatz of the underlying exactly solvable lattice model.
The term fermionic formula was coined by the Stony Brook
group \cite{KKMM:1993,KKMM:1993a}, who interpreted fermionic-type
formulas for characters and branching functions of conformal
field theory models as partition functions of quasiparticle
systems with ``fractional'' statistics obeying Pauli's exclusion
principle. For type $A_n^{(1)}$ the fermionic formulas were proven
in \cite{KSS:2002} using a generalization of a bijection between
crystals and rigged configurations of Kirillov and 
Reshetikhin~\cite{KR:1988}. In~\cite{OSS:2002} similar bijections were
used to prove the fermionic formula for nonexceptional types for crystals
$B^{1,1}$. Rigged configurations are combinatorial objects which
label the solutions to the Bethe equations.
The bijection between crystals and rigged configurations reflects two 
different methods to solve lattice models in statistical mechanics: 
the corner-transfer-matrix method and the Bethe Ansatz.

The theory of virtual crystals \cite{OSS,OSS:2002a} provides a realization 
of crystals of type $X$ as crystals of type $Y$, based on well-known natural 
embeddings $X\hookrightarrow Y$ of affine algebras~\cite{JM:1985}:
\begin{equation*}
\begin{array}{cll}
C_n^{(1)}, A_{2n}^{(2)}, A_{2n}^{(2)\dagger},D_{n+1}^{(2)}
&\hookrightarrow
&A_{2n-1}^{(1)}\\
A_{2n-1}^{(2)}, B_n^{(1)} &\hookrightarrow &D_{n+1}^{(1)}\\
E_6^{(2)}, F_4^{(1)} &\hookrightarrow &E_6^{(1)}\\
D_4^{(3)}, G_2^{(1)} &\hookrightarrow &D_4^{(1)}.
\end{array}
\end{equation*}
Note that under these embeddings every affine Kac--Moody
algebra is embedded into one of simply-laced type $A_n^{(1)}$,
$D_n^{(1)}$ or $E_6^{(1)}$. Hence, by the virtual
crystal method the combinatorial structure of any finite-dimensional
affine crystal can be expressed in terms of the combinatorial
crystal structure of the simply-laced types.
Whereas the affine crystals $B^{r,s}$ of type $A_n^{(1)}$
are already well-understood \cite{Sh:2002}, this is not the
case for $B^{r,s}$ of types $D_n^{(1)}$ and $E_6^{(1)}$.

In this paper we discuss the affine crystals $B^{r,1}$ for
type $D_n^{(1)}$ and the corresponding rigged configurations.
The crystal structure of $B^{r,1}$ of type $D_n^{(1)}$ is given
by Koga \cite{K:1999}. However, here we need a different
combinatorial description of these crystals compatible with
rigged configurations. 
We present a bijection between crystal paths and rigged configurations
of type $D_n^{(1)}$, thereby proving the fermionic formulas 
of~\cite{HKOTY:1999,HKOTT:2002} for type $D_n^{(1)}$ corresponding to 
tensor products of crystals of the form $B^{r,1}$.

The paper is structured as follows. In Section~\ref{sec:crystal}
the main definitions regarding crystals are reviewed. Section~\ref{sec:crystal D}
deals with the finite-dimensional crystals of type $D_n^{(1)}$ explicitly.
In particular $B^{k,1}$ as a set is given by single column tableaux of
height $k$ with the action of the Kashiwara operators as specified in
Theorem~\ref{thm:affine}. The bijection between crystal paths and rigged
configurations is stated in Section~\ref{sec:bij}. The proof of the main
Theorem~\ref{thm:bij} is relegated to Appendix~\ref{app:bij proof}.
In Section~\ref{sec:statistics} statistics are defined on both crystal paths
and rigged configurations. It is shown in Theorem~\ref{thm:stat} that the
bijection preserves the statistics. This yields a proof of the fermionic
formulas as presented in Section \ref{sec:fermionic}. Appendices \ref{app:emb proof}
and \ref{app:delta-delta} are devoted to proofs of various lemmas.
In Appendix \ref{app:R matrix} the combinatorial $R$-matrix is given explicitly
for highest weight elements.

\subsection*{Acknowledgments}
I am grateful to Mark Shimozono for sharing with me the formulation of
$\dual$-dual crystals for general types presented in Section \ref{sec:dual crystal}
(see also \cite{SS}) and to Lipika Deka for her collaboration in
the early stages of this project~\cite{SD:2003}.
I would also like to thank the Max-Planck-Institut f{\"ur} Mathematik in Bonn
where part of this work was carried out.

\section{Crystal bases} \label{sec:crystal}

\subsection{Axiomatic definition of crystals}
Let $\gggg$ be a symmetrizable Kac-Moody algebra, $P$ the weight
lattice, $I$ the index set for the vertices of the Dynkin diagram
of $\gggg$, $\{\alpha_i\in P \mid i\in I \}$ the simple roots, and
$\{h_i\in P^* \mid i\in I \}$ the simple coroots.
Let $U_q(\gggg)$ be the quantized universal enveloping algebra of
$\gggg$ \cite{K:1995}. A $U_q(\gggg)$-crystal is a nonempty set $B$
equipped with maps $\wt:B\rightarrow P$ and
$e_i,f_i:B\rightarrow B\cup\{\es\}$ for all $i\in I$,
satisfying
\begin{align}
\label{eq:e-f}
f_i(b)=b' &\Leftrightarrow e_i(b')=b
\text{ if $b,b'\in B$} \\
\wt(f_i(b))&=\wt(b)-\alpha_i \text{ if $f_i(b)\in B$} \\
\label{eq:string length}
\inner{h_i}{\wt(b)}&=\varphi_i(b)-\epsilon_i(b).
\end{align}
Here for $b \in B$
\begin{equation*}
\begin{split}
\epsilon_i(b)&= \max\{n\ge0\mid e_i^n(b)\not=\es \} \\
\varphi_i(b) &= \max\{n\ge0\mid f_i^n(b)\not=\es \}.
\end{split}
\end{equation*}
(It is assumed that $\varphi_i(b),\epsilon_i(b)<\infty$ for all
$i\in I$ and $b\in B$.) A $U_q(\gggg)$-crystal $B$ can be viewed
as a directed edge-colored graph (the crystal graph) whose
vertices are the elements of $B$, with a directed edge from $b$ to
$b'$ labeled $i\in I$, if and only if $f_i(b)=b'$.

\subsection{Tensor products of crystals}
Let $B_1,B_2,\dotsc,B_L$ be $U_q(\gggg)$-crystals. The
Cartesian product $B_L\times \dotsm \times B_2 \times B_1$ has the
structure of a $U_q(\gggg)$-crystal using the so-called signature
rule. The resulting crystal is denoted
$B=B_L\otimes\dots\otimes B_2\otimes B_1$ and its elements
$(b_L,\dotsc,b_1)$ are written $b_L\otimes \dotsm \otimes b_1$
where $b_j\in B_j$. The reader is warned that our convention is
opposite to that of Kashiwara \cite{K:1995}. Fix $i\in I$ and
$b=b_L\otimes\dotsm\otimes b_1\in B$. The $i$-signature of $b$ is
the word consisting of the symbols $+$ and $-$ given by
\begin{equation*}
\underset{\text{$\varphi_i(b_L)$ times}}{\underbrace{-\dotsm-}}
\quad \underset{\text{$\epsilon_i(b_L)$
times}}{\underbrace{+\dotsm+}} \,\dotsm\,
\underset{\text{$\varphi_i(b_1)$ times}}{\underbrace{-\dotsm-}}
\quad \underset{\text{$\epsilon_i(b_1)$
times}}{\underbrace{+\dotsm+}} .
\end{equation*}
The reduced $i$-signature of $b$ is the subword of the
$i$-signature of $b$, given by the repeated removal of adjacent
symbols $+-$ (in that order); it has the form
\begin{equation*}
\underset{\text{$\varphi$ times}}{\underbrace{-\dotsm-}} \quad
\underset{\text{$\epsilon$ times}}{\underbrace{+\dotsm+}}.
\end{equation*}
If $\varphi=0$ then $f_i(b)=\es$; otherwise
\begin{equation*}
f_i(b_L\otimes\dotsm\otimes b_1)= b_L\otimes \dotsm \otimes
b_{j+1} \otimes f_i(b_j)\otimes \dots \otimes b_1
\end{equation*}
where the rightmost symbol $-$ in the reduced $i$-signature of
$b$ comes from $b_j$. Similarly, if $\epsilon=0$ then
$e_i(b)=\es$; otherwise
\begin{equation*}
e_i(b_L\otimes\dotsm\otimes b_1)= b_L\otimes \dotsm \otimes
b_{j+1} \otimes e_i(b_j)\otimes \dots \otimes b_1
\end{equation*}
where the leftmost symbol $+$ in the reduced $i$-signature of $b$
comes from $b_j$. It is not hard to verify that this well-defines
the structure of a $U_q(\gggg)$-crystal with
$\varphi_i(b)=\varphi$ and $\epsilon_i(b)=\epsilon$ in the above
notation, with weight function
\begin{equation} \label{eq:tensor wt}
\wt(b_L\otimes\dotsm\otimes b_1)=\sum_{j=1}^L \wt(b_j).
\end{equation}
This tensor construction is easily seen to be associative. The
case of two tensor factors is given explicitly by
\begin{equation} \label{eq:f on two factors}
f_i(b_2\otimes b_1) = \begin{cases} f_i(b_2)\otimes b_1
& \text{if $\epsilon_i(b_2)\ge \varphi_i(b_1)$} \\
b_2\otimes f_i(b_1) & \text{if $\epsilon_i(b_2)<\varphi_i(b_1)$}
\end{cases}
\end{equation}
and
\begin{equation} \label{eq:e on two factors}
e_i(b_2\otimes b_1) = \begin{cases} e_i(b_2) \otimes b_1 &
\text{if $\epsilon_i(b_2)>\varphi_i(b_1)$} \\
b_2\otimes e_i(b_1) & \text{if $\epsilon_i(b_2)\le \varphi_i(b_1)$.}
\end{cases}
\end{equation}

\subsection{Combinatorial $R$-matrix}\label{sec:R}
Let $B_1$ and $B_2$ be two $U'_q(\gggg)$-crystals associated with 
finite-dimensional $U'_q(\gggg)$-modules such that $B_1\otimes B_2$
and $B_2\otimes B_1$ are connected and such that there is a weight $\la\in P$ with 
unique elements $u\in B_1\otimes B_2$ and $\tilde{u}\in B_2\otimes B_1$ of weight 
$\wt(u)=\wt(\tilde{u})=\la$. Then there is a unique crystal isomorphism
\begin{equation*}
R:B_1\otimes B_2\to B_2\otimes B_1,
\end{equation*}
called the combinatorial $R$-matrix. Explicitly, $R(u)=\tilde{u}$ and the action of 
$R$ on any other element in $B_1\otimes B_2$ is determined by using
\begin{equation}
\label{eq:R crystal iso}
 R(e_i(b))=e_i(R(b)) \qquad \text{and} \qquad R(f_i(b))=f_i(R(b)).
\end{equation}

\section{The crystal $B^{k,1}$ of type $D_n^{(1)}$} \label{sec:crystal D}

{}From now on we restrict our attention to crystals of type $D_n$ and
$D_n^{(1)}$.

\subsection{Dynkin data}
For type $D_n$, the simple roots are
\begin{equation*}
\begin{split}
\alpha_i&=\varepsilon_i-\varepsilon_{i+1} \qquad \text{for $1\le i<n$}\\
\alpha_n&=\varepsilon_{n-1}+\varepsilon_n
\end{split}
\end{equation*}
and the fundamental weights are
\begin{equation*}
\begin{aligned}
\La_i&=\varepsilon_1+\cdots+\varepsilon_i &\text{for $1\le i\le n-2$}\\
\La_{n-1}&=(\varepsilon_1+\cdots+\varepsilon_{n-1}-\varepsilon_n)/2 &\\
\La_n&=(\varepsilon_1+\cdots+\varepsilon_{n-1}+\varepsilon_n)/2&
\end{aligned}
\end{equation*}
where $\varepsilon_i\in Z^n$ is the $i$-th unit standard vector.
We also define
\begin{equation*}
\begin{aligned}
\om_i&=\La_i=\varepsilon_1+\cdots+\varepsilon_i & \text{for $1\le i\le n-2$}\\
\om_{n-1}&=\La_{n-1}+\La_n=\varepsilon_1+\cdots+\varepsilon_{n-1}&\\
\om_n&=2\La_n=\varepsilon_1+\cdots+\varepsilon_n&\\
\omb_n&=2\La_{n-1}=\varepsilon_1+\cdots+\varepsilon_{n-1}-\varepsilon_n.&
\end{aligned}
\end{equation*}

\subsection{The crystals of antisymmetric tensor representations} \label{sec:column}
\begin{figure}
\scalebox{0.9}{
\begin{picture}(365,100)(-10,-50)
\BText(0,0){1} \LongArrow(10,0)(30,0) \BText(40,0){2}
\LongArrow(50,0)(70,0) \Text(85,0)[]{$\cdots$}
\LongArrow(95,0)(115,0) \BText(130,0){n-1}
\LongArrow(143,2)(160,14) \LongArrow(143,-2)(160,-14)
\BText(170,15){n} \BBoxc(170,-15)(13,13)
\Text(170,-15)[]{\footnotesize$\overline{\mbox{n}}$}
\LongArrow(180,14)(197,2) \LongArrow(180,-14)(197,-2)
\BBoxc(215,0)(25,13)
\Text(215,0)[]{\footnotesize$\overline{\mbox{n-1}}$}
\LongArrow(230,0)(250,0) \Text(265,0)[]{$\cdots$}
\LongArrow(275,0)(295,0) \BBoxc(305,0)(13,13)
\Text(305,0)[]{\footnotesize$\overline{\mbox{2}}$}
\LongArrow(315,0)(335,0) \BBoxc(345,0)(13,13)
\Text(345,0)[]{\footnotesize$\overline{\mbox{1}}$}
\PText(20,2)(0)[b]{1}
\PText(60,2)(0)[b]{2} \PText(105,2)(0)[b]{n-2}
\PText(152,13)(0)[br]{n-1} \PText(152,-9)(0)[tr]{n}
\PText(188,13)(0)[bl]{n} \PText(188,-9)(0)[tl]{n-1}
\PText(240,2)(0)[b]{n-2} \PText(285,2)(0)[b]{2}
\PText(325,2)(0)[b]{1} 
\end{picture}
}
\caption{\label{fig:vr}Crystal $B(\om_1)$ of the vector representation}
\end{figure}
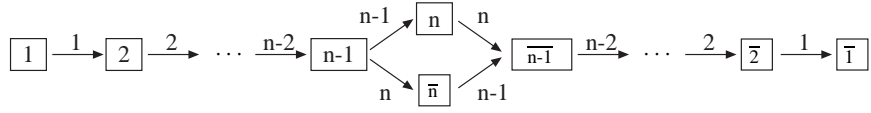

The $U_q(D_n)$-crystal graphs of the antisymmetric tensor representations \cite{KN:1994}
are denoted by $B(\om_\ell)$ (resp. $B(\omb_n)$).
The crystal graph $B(\om_1)$ of the vector representation is given
in Figure \ref{fig:vr}.
The crystal $B(\om_\ell)$ (resp. $B(\omb_n)$) is the connected component
of $B(\om_1)^{\otimes \ell}$ (resp. $B(\om_1)^{\otimes n}$) containing
the element $\ell\otimes \ell-1\otimes \cdots \otimes 1$ 
(resp. $\overline{n}\otimes n-1\otimes \cdots \otimes 1$).

Explicitly, the elements of $B(\om_\ell)$ for $1\le \ell<n$
are single columns of height $\ell$ with entries from the alphabet
$\{1<2<\cdots<n-1<\substack{n\\ \overline{n}}<\overline{n-1}<\cdots 
<\overline{1}\}$. In Young tableau notation a column is represented by
\begin{equation*}
\begin{array}{|c|}
\hline m_1 \\ \hline m_2\\ \hline \vdots \\ \hline m_\ell\\ \hline \end{array}
\qquad \text{which we abbreviate by $m_\ell\ldots m_1$.}
\end{equation*}
If $m_\ell\ldots m_1$ is in $B(\om_\ell)$ it has to satisfy
\begin{equation}\label{eq:sc}
\begin{split}
&\text{1. $(m_j,m_{j+1})=(n,\overline{n}),(\overline{n},n)$
or $m_j<m_{j+1}$ for all $j$;}\\
&\text{2. If $m_a=p$ and $m_b=\overline{p}$, then $\dist(\overline{p},p)\le p$.}
\end{split}
\end{equation}
Here $\dist(\overline{p},p)=\ell+1+a-b$ if $m_a=p$ and $m_b=\overline{p}$.

The elements of $B(\om_n)$ (resp. $B(\omb_n)$) are single columns of
height $n$. An element $m_n\ldots m_1\in B(\om_n)$ (resp. $B(\omb_n)$)
satisfies \eqref{eq:sc} and in addition
\begin{equation}\label{eq:n}
\begin{split}
&\text{3. If $m_j=n$, then $n-j$ is even (resp. odd);}\\
&\text{\phantom{3.} if $m_j=\overline{n}$, then $n-j$ is odd (resp. even).}
\end{split}
\end{equation}

The action of $e_i,f_i$ for $1\le i\le n$ is given by embedding
$B(\om_\ell)\hookrightarrow B(\om_1)^{\otimes \ell}$ (resp.
$B(\omb_n)\hookrightarrow B(\om_1)^{\otimes n}$) in the natural way
by mapping $m_\ell\cdots m_1$ to $m_\ell\otimes \cdots \otimes m_1$
and using the tensor product rules.

\subsection{The crystals of spinor representations}
There are two crystals $B(\La_n)$ and $B(\La_{n-1})$ associated
with the spinor representations of $D_n$ \cite{KN:1994}.
The elements of $B(\La_n)$ (resp. $B(\La_{n-1})$) are given by
single columns $m_n\ldots m_1$ such that
\begin{equation}\label{eq:spinor}
\begin{split}
&\text{1. $m_j<m_{j+1}$;}\\
&\text{2. $m$ and $\overline{m}$ cannot occur simultaneously;}\\
&\text{3. If $m_j=n$, then $n-j$ is even (resp. odd);}\\
&\text{\phantom{3.} if $m_j=\overline{n}$, then $n-j$ is odd (resp. even).}
\end{split}
\end{equation}
If $b\in B(\La_n)$ or $B(\La_{n-1})$ and $b$ contains $i$ and 
$\overline{i+1}$ for $1\le i<n$, then $f_ib$ is obtained by replacing 
$i$ by $i+1$ and $\overline{i+1}$ by $\overline{i}$. Else $f_ib=0$.
If $b\in B(\La_n)$ or $B(\La_{n-1})$ and $b$ contains $n-1$ and $n$, 
then $f_n b$ is obtained by replacing $n-1$ by $\overline{n}$ and 
$n$ by $\overline{n-1}$. Otherwise $f_n b=0$.

Note that for $k=n-1,n$ there is the following embedding of crystals
\begin{equation*}
B(\La_k) \hookrightarrow B(2\La_k)
\end{equation*}
where $e_i$ and $f_i$ act as $e_i^2$ and $f_i^2$, respectively.

\subsection{Affine crystals}
As a classical crystal the affine crystal $B^{k,1}$ is isomorphic to
\begin{equation}\label{eq:classical comp}
B^{k,1}\cong \begin{cases} 
B(\Lambda_k)\oplus B(\Lambda_{k-2}) \oplus \cdots \oplus B(0) 
 & \text{if $k$ is even, $1\le k\le n-2$}\\
B(\Lambda_k)\oplus B(\Lambda_{k-2}) \oplus \cdots \oplus B(\Lambda_1) 
 & \text{if $k$ is odd, $1\le k\le n-2$}\\
B(\Lambda_k) & \text{if $k=n-1,n$.}
\end{cases}
\end{equation}

The rule for $f_0$ on $B^{k,1}$ for $k=n,n-1$ is as follows.
If $b=\overline{1}\overline{2}m_{n-2}\ldots m_1$ then
$f_0 b=m_{n-2}\ldots m_1 21$. Otherwise $f_0b=0$.
The operator $f_0$ on $B^{k,1}$ when $1\le k<n-1$ was given
by Koga \cite{K:1999} in terms of a tensor product of two spinor
representations. For our purposes, we will need a combinatorial
description of $f_0$ on elements $b\in B^{k,1}$ represented
by a column of height $k$. This description of $b$ will also be
used in the description of the bijection from crystals to rigged
configurations.

For $1\le k\le n-2$, we want to represent $B^{k,1}$ as the set of all column
tableaux of height $k$ satisfying \eqref{eq:sc} point 1. By \eqref{eq:classical comp}
$B^{k,1}$ is the union of the sets corresponding to $B(\La_{\ell})$ with
$k-\ell \equiv 0 \pmod{2}$. We will define maps from $B(\La_\ell)$ to
column tableaux of height $k$.
As explained in Section \ref{sec:column}, $b\in B(\La_\ell)$ can be 
represented by a column of height $\ell$. If $b\in B(\La_\ell)\subset B^{k,1}$
then fill the column of height $\ell$ of $b$ successively by a pair 
$(\overline{i}_j,i_j)$ for $1\le j\le (k-\ell)/2$ in the following way to 
obtain a column of height $k$. Set $i_0=0$.
Let $i_{j-1}<i_j$ be minimal such that 
\begin{enumerate}
\item
neither $i_j$ nor $\overline{i_j}$ is in the column;
\item
adding $i_j$ and $\overline{i_j}$ to the column we have
$\dist(\overline{i_j},i_j)\ge i_j+j$;
\item
adding $i_j$ and $\overline{i_j}$ to the column, all other pairs 
$(\overline{a},a)$ in the new column with $a>i_j$ 
satisfy $\dist(\overline{a},a)\le a+j$.
\end{enumerate}
It is not hard to see that $f_i$ and the filling map commute, that is, the
classical crystal graph with edges $f_1,\ldots,f_n$ does not change.
Denote the filling map to height $k$ by $\F_k$ or simply $\F$.

Let $\D_k$ or $\D$, the dropping map, be the inverse of $\F_k$. 
Explicitly, given a one-column tableau $b$ of height $k$, let $i_0=0$ and
successively find $i_j>i_{j-1}$ minimal such that the pair $(\overline{i_j},i_j)$
is in $b$ and $\dist(\overline{i_j},i_j)\ge i_j+j$. Drop all such pairs
$(\overline{i_j},i_j)$ from $b$.

For the definition of $f_0$ we need slight variants of $\F_k$ and $\D_k$, 
denoted by $\Ft_k$ and $\Dt_k$, respectively, which act on columns that do not contain
$1,2,\overline{2},\overline{1}$. On these columns $\Ft_k$ and $\Dt_k$ are defined 
by replacing $i\mapsto i-2$ and $\overline{i}\mapsto \overline{i-2}$, then
applying $\F_k$ and $\D_k$, and finally replacing $i\mapsto i+2$ and 
$\overline{i}\mapsto \overline{i+2}$.

\begin{example}
Let $b=\overline{3}\overline{5}532\in B(\Lambda_5)\subset B^{9,1}$.
Then $\F_9(b)=\overline{3567}76532$. Similarly if 
$b=\overline{2345}65421\in B^{9,1}$ then $\D_9(b)=\overline{3}61$.
Finally, $\Ft_5(\overline{4}64)=\overline{45}654$.
\end{example}

\begin{theorem}\label{thm:affine}
Let $b\in B^{k,1}$. Then
\begin{equation}
f_0b=\begin{cases}
 \F_k(\Dt_{k-2}(x)) & \text{if $b=\overline{1}\overline{2}x$}\\
 \Ft_{k-1}(x)2 & \text{if $b=\overline{1}\overline{2}x2$}\\
 \Ft_{k-1}(x)1 & \text{if $b=\overline{1}\overline{2}x1$}\\
 \Ft_{k-2}(x)21 & \text{if $b=\overline{1}\overline{2}x21$}\\
 \F_k(\Dt_{k-1}(x)2) & \text{if $b=\overline{1}x$}\\
 \F_k(\Dt_{k-1}(x)1) & \text{if $b=\overline{2}x$}\\
 x21 & \text{if $b=\overline{1}x1$ and $\Dt_{k-2}(x)=x$}\\
 0 & \text{otherwise}
\end{cases}
\end{equation}
\begin{equation}
e_0b=\begin{cases}
 \F_k(\Dt_{k-2}(x)) & \text{if $b=x21$}\\
 \overline{2}\Ft_{k-1}(x) & \text{if $b=\overline{2}x21$}\\
 \overline{1}\Ft_{k-1}(x) & \text{if $b=\overline{1}x21$}\\
 \overline{1}\overline{2}\Ft_{k-2}(x) & \text{if $b=\overline{1}\overline{2}x21$}\\
 \F_k(\overline{2}\Dt_{k-1}(x)) & \text{if $b=x1$}\\
 \F_k(\overline{1}\Dt_{k-1}(x)) & \text{if $b=x2$}\\
 \overline{1}\overline{2}x & \text{if $b=\overline{1}x1$ and $\Dt_{k-2}(x)=x$}\\
 0 & \text{otherwise}
\end{cases}
\end{equation}
where $x$ does not contain $1,2,\overline{2},\overline{1}$.
\end{theorem}

\begin{proof}
In \cite{K:1999}, Koga gives the crystal $B^{k,1}$ in terms of a 
tensor product of two spinor crystals.
Let us denote the map from our description in terms of columns
of height $k$ to Koga's picture by $\psi$. Since $F_k$ commutes
with all $e_i,f_i$ for $1\le i\le n$, $B^{k,1}\cong \psi(B^{k,1})$
as classical crystals. Hence $\psi$ is specified by its action
on highest weight vectors which is given by
\begin{equation*}
\psi(\overline{\ell+1}\;\cdots\overline{p}\; p\cdots \ell+1\; \ell\cdots 21)
=\overline{\ell+1}\cdots \overline{n}\;\ell\cdots 21\; \otimes\; n\cdots 21
\end{equation*}
where $p=(k+\ell)/2$.

Note that $f_0$ commutes with $f_i$ for $3\le i\le n$. Hence,
to show that the above action of $f_0$ on $B^{k,1}$ is equivalent
to the action given by Koga, it suffices to show that
$\psi(f_0 b)=f_0 \psi(b)$ where $b$ is highest weight with respect
to $e_i$ for $3\le i\le n$, that is $e_ib=0$ for these $i$. Call
such $b$ $J'$-highest weight where $J'=\{3,4,\ldots,n\}$. This can
be done explicitly for each case.

\end{proof}

Alternatively, the action of $f_0$ and $e_0$ is given by 
$f_0=\sigma\circ f_1\circ \sigma$ and $e_0=\sigma\circ e_1\circ \sigma$
where $\sigma$ is an involution on $B^{k,1}$ based on the type $D_n^{(1)}$ 
Dynkin diagram automorphism that interchanges nodes 0 and 1. If $b\in B^{k,1}$
contains either 1 or $\overline{1}$, $\sigma(b)$ is obtained from
$b$ by interchanging 1 and $\overline{1}$. If $b=\overline{1}x1$,
then $\sigma(b)=\Fh_k(\D(x))$. If $b$ contains neither 1 nor
$\overline{1}$ then $\sigma(b)=F_k(\Dh(b))$.
Here $\Dh$ and $\Fh$ act on columns not containing $1$ and $\overline{1}$ and are 
defined in the same way as $\Dt$ and $\Ft$
with the maps $i\mapsto i\pm 2$, $\overline{i}\mapsto \overline{i\pm 2}$
replaced by $i\mapsto i\pm 1$, $\overline{i}\mapsto \overline{i\pm 1}$.

For later use, it will be convenient to define the sets $\Bh^{k,1}$ for $1\le k\le n$
and $\Bbh^{n,1}$. For $1\le k\le n-2$, $\Bh^{k,1}=B^{k,1}$. $\Bh^{n-1,1}$ is the set
of all single columns of height $n-1$ satisfying condition \eqref{eq:sc} point 1.
$\Bh^{n,1}$ (resp. $\Bbh^{n,1}$) is the set of all single columns of height $n$
satisfying \eqref{eq:sc} point 1 and \eqref{eq:n}. An affine crystal structure
can be defined on $\Bh^{n-1,1}$, $\Bh^{n,1}$ and $\Bbh^{n,1}$ by extending
the rules of Theorem \ref{thm:affine}.
\begin{theorem}
As affine crystals
\begin{equation*}
\begin{split}
 \Bh^{n-1,1} &\cong B^{n,1}\otimes B^{n-1,1}\\
 \Bh^{n,1} &\cong B^{n,1}\otimes B^{n,1}\\
 \Bbh^{n,1} &\cong B^{n-1,1}\otimes B^{n-1,1}.
\end{split}
\end{equation*}
\end{theorem}
\begin{proof} The proof is similar to the proof of Theorem \ref{thm:affine}.
\end{proof}

There are affine crystal embeddings
\begin{equation}\label{eq:emb affine}
 \emb_B: B^{k,1} \hookrightarrow B^{k,1} \otimes B^{k,1}
\end{equation}
mapping $u(B^{k,1})\mapsto u(B^{k,1})\otimes u(B^{k,1})$ for $1\le k\le n$ and sending 
$e_i$ and $f_i$ to $e_i^2$ and $f_i^2$, respectively, for $0\le i\le n$. Here 
$u(B^{k,1})=k(k-1)\cdots 1$ for $1\le k\le n-2$, $u(B^{n-1,1})=\overline{n}(n-1)\cdots 21$ and 
$u(B^{n,1})=n(n-1)\cdots 21$. Note that the elements in the image of this embedding are aligned
(see \cite[Section 6.4]{OSS}), meaning in this case that $\varphi_i(b)$ and
$\epsilon_i(b)$ are even for all $0\le i\le n$. Denote the image of $B^{k,1}$ under
this embedding by $E^{k,1}=\emb_B(B^{k,1})$. For $k=n-1,n$ we view 
$E^{n-1,1}\subset \Bbh^{n,1}$ and $E^{n,1}\subset \Bh^{n,1}$.

As classical crystals
\begin{equation*}
\begin{split}
 \Bh^{n-1,1} &\cong B(\om_{n-1}) \oplus B(\om_{n-3}) \oplus B(\om_{n-5}) \oplus \cdots\\
 \Bh^{n,1} &\cong B(\om_n) \oplus B(\om_{n-2}) \oplus B(\om_{n-4}) \oplus \cdots\\
 \Bbh^{n,1} &\cong B(\omb_n) \oplus B(\om_{n-2}) \oplus B(\om_{n-4}) \oplus \cdots.
\end{split} 
\end{equation*}

\subsection{Paths}
Let $B=B^{k_L,1}\otimes B^{k_{L-1},1}\otimes \cdots \otimes B^{k_1,1}$ be
a type $D_n^{(1)}$ crystal with $1\le k_i\le n$. For a dominant
integral weight $\la$ define the set of paths as follows
\begin{equation*}
\Path(\la,B)=\bigl\{ b\in B \mid 
 \text{$\wt(b)=\la$, $e_i(b)=\es$ for all $i\in J$} \bigr\}
\end{equation*}
where $J=I\backslash \{0\}=\{1,2,\ldots,n\}$.

\subsection{Dual crystals}\label{sec:dual crystal}
Let $\omega_0$ be the longest element in the Weyl group of $D_n$.
The action of $\omega_0$ on the weight lattice $P$ of $D_n$ is given by
\begin{equation*}
\begin{split}
\omega_0(\La_i) &= -\La_{\tau(i)}\\
\omega_0(\alpha_i) &= -\alpha_{\tau(i)}
\end{split}
\end{equation*}
where $\tau:J\to J$ is the identity if $n$ is even and interchanges
$n-1$ and $n$ and fixes all other Dynkin nodes if $n$ is odd. The automorphism
$\tau$ extends to the Dynkin diagram of $D_n^{(1)}$ by defining
$\tau(0)=0$.

There is a unique involution $\dual:B\to B$, called dual map, satisfying
\begin{equation*}
\begin{split}
\wt(b^\dual) &= \omega_0 \wt(b)\\
e_i(b)^\dual &= f_{\tau(i)}(b^\dual)\\
f_i(b)^\dual &= e_{\tau(i)}(b^\dual).
\end{split}
\end{equation*}
Let $u(B)$ be the highest weight element in the highest component of $B$
(for $B=B^{k,1}$ with $1\le k\le n-2$ the highest component is $B(\La_k)$).
Uniqueness of $\dual$ follows from the connectedness of $B$ and the fact that 
$u(B)$ is the unique vector of its weight. We have 
\begin{equation*}
(B_1\otimes B_2)^\dual \cong B_2\otimes B_1
\end{equation*}
with $(b_1\otimes b_2)^\dual \mapsto b_2^\dual \otimes b_1^\dual$.

Explicitly, on $B(\La_1)$ the involution $\dual$ is given by 
\begin{equation*}
 i \longleftrightarrow \overline{i}
\end{equation*}
except for $i=n$ with $n$ odd in which case $n \leftrightarrow n$ and 
$\overline{n} \leftrightarrow \overline{n}$. A column $m_k\cdots m_1\in B^{k,1}$
is mapped to $m_1^\dual\cdots m_k^\dual$.

\section{Bijection between paths and rigged configurations} \label{sec:bij}

Let $B=B^{k_L,1}\otimes B^{k_{L-1},1}\otimes \cdots \otimes B^{k_1,1}$ be
a type $D_n^{(1)}$ crystal and $\la$ a dominant integral weight.
In this section we define the set of rigged configurations
$\RC(\la,B)$ and a bijection $\Phi:\RC(\la,B)\to \Path(\la,B)$.

\subsection{Rigged configurations} \label{sec:rc}
Denote by $L^a$ the number of factors $B^{a,1}$ in $B$,
so that $B=\bigotimes_{a=1}^n (B^{a,1})^{\otimes L^a}$.
For later use, let us also define $\Lh^{n-1}$ (resp. $\Lh^n$, $\Lbh^n$) to be the 
number of tensor factors $\Bh^{n-1,1}$ (resp. $E^{n,1}$, $E^{n-1,1}$) in $B$.
Let
\begin{equation}\label{eq:L}
 \vL = \sum_{a\in J} L^a \La_a + \sum_{a\in\{n-1,n\}} \Lh^a \om_a + \Lbh^n \omb_n.
\end{equation}
Then the sequence of partitions $\nt=(\nu^{(1)},\nu^{(2)},\ldots,\nu^{(n)})$
is a $(\la,B)$-configuration if
\begin{equation}\label{eq:size}
\sum_{\substack{a\in J\\ i\ge 1}} i m_i^{(a)} \alpha_a
= \vL - \la,
\end{equation}
where $m_i^{(a)}$ is the number of parts of length $i$ in 
partition $\nu^{(a)}$. 
%
%
A $(\la,B)$-configuration is admissible if $p_i^{(a)}\ge 0$ for all 
$a\in J$ and $i\ge 1$, where $p_i^{(a)}$ is the vacancy number
\begin{equation}\label{eq:vac}
p_i^{(a)}=\bigl(\alpha_a \mid \vL 
 - \sum_{b\in J} \alpha_b \sum_{j\ge 1} \min(i,j)m_j^{(b)}\bigr),
\end{equation}
where $(\cdot \mid \cdot )$ is the normalized invariant form on $P$
such that $(\alpha_i\mid \alpha_j)$ is the Cartan matrix.
%
%
Denote the set of admissible $(\la,B)$-configurations by
$\Conf(\la,B)$. A rigged configuration $(\nt,\Jt)$ is a sequence
of partitions $\nt=\{\nu^{(a)}\}_{a\in J} \in \Conf(\la,B)$ together with a
double sequence of partitions $\Jt=\{J^{(a,i)}\}_{\substack{a\in J\\
i\ge 1}}$ such that the partition $J^{(a,i)}$ is contained
in a box of size $m_i^{(a)}\times p_i^{(a)}$. The set of rigged
configurations is denoted by $\RC(\la,B)$.

The partition $J^{(a,i)}$ is called singular if it has a part of size
$p_i^{(a)}$. A singular string is a pair $(i,p)$ where $i$ is the part
of the partition $\nu^{(a)}$ and $p=p_i^{(a)}$ is the size of the largest
part of $J^{(a,i)}$.

\subsection{The bijection $\Phi$}
Let $B'$ be a tensor product of type $D_n^{(1)}$ crystals of the 
form $B^{a,1}$ or $\Bh^{n-1,1}$.
\begin{theorem}\label{thm:bij}
There exists a unique family of bijections $\Phi_B:\RC(\la,B)\to\Path(\la,B)$ 
such that:
\begin{enumerate}
\item \label{pt:single}
Suppose $B=B^{1,1}\otimes B'$. Then the diagram
\begin{equation*}
\begin{CD}
\RC(\la, B)        @>{\Phi_B}>> \Path(\la,B)\\
@V{\delta}VV @VV{\lh}V\\
\bigcup_{\la^-} \RC(\la^-,B')  @>{\Phi_{B'}}>>\bigcup_{\la^-}\Path(\la^-,B')
\end{CD}
\end{equation*}
commutes.
\item \label{pt:split}
Suppose $B=\Bh^{k,1}\otimes B'$ with $2\le k\le n-1$ 
and let $\Bt=B^{1,1}\otimes \Bh^{k-1,1}\otimes B'$. Then the diagram
\begin{equation*}
\begin{CD}
\RC(\la,B)        @>{\Phi_B}>> \Path(\la,B)\\
@V{\tj}VV @VV{\ts}V\\
\RC(\la,\Bt)  @>{\Phi_{\Bt}}>>\Path(\la,\Bt)
\end{CD}
\end{equation*}
commutes.
\item \label{pt:spinor}
Suppose $B=B^{k,1}\otimes B'$ with $k=n-1,n$. Then the diagram
\begin{equation*}
\begin{CD}
\RC(\la,B)        @>{\Phi_B}>> \Path(\la,B)\\
@V{\ds}VV @VV{\lhs}V\\
\bigcup_{\la^{\dspin}}\RC(\la^{\dspin},B')  @>{\Phi_{B'}}>>
 \bigcup_{\la^{\dspin}}\Path(\la^{\dspin},B')
\end{CD}
\end{equation*}
commutes.
\end{enumerate} 
\end{theorem}
A proof of this theorem is given in Appendix \ref{app:bij proof}.
In the following we define the various yet undefined maps.

\subsection{Various maps}\label{sec:var maps}
The map $\lh:\Path(\la,B^{1,1}\otimes B')\to\bigcup_{\la^-}\Path(\la^-,B')$ 
in point \ref{pt:single} of Theorem \ref{thm:bij} is given by 
$b\otimes b'\mapsto b'$. The union is over all
$\la^-$ such that $\la-\la^-$ is the weight of an element in $B^{1,1}$.

Similarly, the map $\lhs:\Path(\la,B)\to\bigcup_{\la^{\dspin}}
\Path(\la^{\dspin},B')$ of point \ref{pt:spinor} is given by
$b\otimes b'\mapsto b'$. The union is over all
$\la^{\dspin}$ such that $\la-\la^{\dspin}$ is the weight of an 
element in $B^{k,1}$ for $k=n-1$ or $n$.

The map top split $\ts:\Path(\la,B) \to \Path(\la,\Bt)$ in point \ref{pt:split} 
of the theorem is obtained by sending $m_k m_{k-1} \cdots m_1 \otimes b'
\in \Path(\la,\Bh^{k,1}\otimes B')$ to $m_k \otimes  m_{k-1} \cdots m_1 
\otimes b' \in \Path(\la,B^{1,1}\otimes \Bh^{k-1,1}\otimes B')$ for 
$1\le k\le n-1$. Similarly we may also define $\ts$ to map
$m_n m_{n-1}\cdots m_1\otimes b'\in \Path(\la,E^{n-1,1}\otimes B')$ 
or $\Path(\la,E^{n,1}\otimes B')$ to
$m_n\otimes m_{n-1}\cdots m_1\otimes b'\in \Path(\la,B^{1,1}\otimes \Bh^{n-1,1}
\otimes B')$.

Define $\tj:\RC(\la,B) \to \RC(\la,\Bt)$ in the following way.
Let $(\nt,\Jt)\in \RC(\la,B)$.
If $B=\Bh^{k,1}\otimes B'$ for $1\le k\le n-1$, $\tj(\nt,\Jt)$ is obtained 
from $(\nt,\Jt)$ by adding a singular string of length 1 to each of 
the first $k-1$ rigged partitions of $\nt$. 
For $B=E^{n,1}\otimes B'$, add a singular
string of length 1 to $(\nt,\Jt)^{(a)}$ for $1\le a\le n-1$.
For $B=E^{n-1,1}\otimes B'$, add a singular
string of length 1 to $(\nt,\Jt)^{(n)}$ and $(\nt,\Jt)^{(a)}$ for $1\le a\le n-2$.

\begin{lemma} \label{lem:tj}
$\tj$ is a well-defined injection that preserves the vacancy numbers.
\end{lemma}
\begin{proof}
To show that $\tj$ is well-defined, we need to check that
$(\ntt,\Jtt)=\tj(\nt,\Jt)\in\RC(\la,\Bt)$ for $(\nt,\Jt)\in \RC(\la,B)$.
Let us first consider the case $B=\Bh^{k,1}\otimes B'$. In this case
$\mt_i^{(a)}=m_i^{(a)}+\chi(1\le a<k)\delta_{i,1}$ and
$\tilde{\vL}=\vL+\om_1+\om_{k-1}-\om_k$. Hence the change of both the left and
right hand side of \eqref{eq:size} is $\sum_{1\le a\le k-1}\alpha_a
=\om_1+\om_{k-1}-\om_k=\epsilon_1-\epsilon_k$.
Therefore \eqref{eq:size} holds for $(\ntt,\Jtt)$.

Using \eqref{eq:vac} the vacancy numbers transform as
\begin{equation*}
\tilde{p}_i^{(a)}=p_i^{(a)}+\bigl(\alpha_a \mid \om_1+\om_{k-1}-\om_k
 - \sum_{1\le b<k} \alpha_b\bigr)=p_i^{(a)},
\end{equation*}
which shows that they are preserved under $\tj$.

The cases $B=E^{n,1}\otimes B'$ and $B=E^{n-1,1}\otimes B'$ are analogous.
\end{proof}

The maps $\delta$ and $\ds$ are defined in the next two subsections.

\subsection{Algorithm for $\delta$}\label{sec:delta}
Let $(\nt,\Jt)\in \RC(\la,B^{1,1}\otimes B')$.
The map $\delta:\RC(\la,B^{1,1}\otimes B')\to \bigcup_{\la^-}
\RC(\la^-,B')$ is defined by the following algorithm \cite{OSS:2002}.
Recall that the partition $J^{(a,i)}$ is called singular if it has a part of size
$p_i^{(a)}$.

Set $\ell^{(0)}=1$ and repeat the following process for $a=1,2,\ldots,n-2$
or until stopped. Find the smallest index $i\ge \ell^{(a-1)}$ such
that $J^{(a,i)}$ is singular. If no such $i$ exists, set $b=a$ and stop.
Otherwise set $\ell^{(a)}=i$ and continue with $a+1$.

If the process has not stopped at $a=n-2$,
find the minimal indices $i,j\ge \ell^{(n-2)}$ such that
$J^{(n-1,i)}$ and $J^{(n,j)}$ are singular. If neither $i$ nor
$j$ exist, set $b=n-1$ and stop.
If $i$ exists, but not $j$, set $\ell^{(n-1)}=i$, $b=n$ and stop.
If $j$ exists, but not $i$, set $\ell^{(n)}=j$, $b=\overline{n}$
and stop. If both $i$ and $j$ exist, set $\ell^{(n-1)}=i$, $\ell^{(n)}=j$
and continue with $a=n-2$.

Now continue for $a=n-2,n-3,\ldots,1$ or until stopped.
Find the minimal index $i\ge \lb^{(a+1)}$ where $\lb^{(n-1)}
=\max(\ell^{(n-1)},\ell^{(n)})$ such that $J^{(a,i)}$ is singular
(if $i=\ell^{(a)}$ then there need to be two parts of size
$p_i^{(a)}$ in $J^{(a,i)}$).
If no such $i$ exists, set $b=\overline{a+1}$ and stop.
If the process did not stop, set $b=\overline{1}$.
Set all yet undefined $\ell^{(a)}$ and $\lb^{(a)}$ to $\infty$.

The rigged configuration $(\ntt,\Jtt)=\delta(\nt,\Jt)$ is obtained by 
removing a box from the selected strings and making the new strings singular
again.
Explicitly
\begin{equation*}
 m_i^{(a)}(\ntt)=m_i^{(a)}(\nt)+\begin{cases}
 1 & \text{if $i=\ell^{(a)}-1$}\\
 -1 & \text{if $i=\ell^{(a)}$}\\
 1 & \text{if $i=\lb^{(a)}-1$ and $1\le a\le n-2$}\\
 -1 & \text{if $i=\lb^{(a)}$ and $1\le a \le n-2$}\\
 0 & \text{otherwise.} \end{cases}
\end{equation*}
If two conditions hold, then all changes should be performed.
The partition $\tilde{J}^{(a,i)}$ is obtained from $J^{(a,i)}$ by removing
a part of size $p_i^{(a)}(\nt)$ for $i=\ell^{(a)}$ and $i=\lb^{(a)}$,
adding a part of size $p_i^{(a)}(\ntt)$ for $i=\ell^{(a)}-1$ and
$i=\lb^{(a)}-1$, and leaving it unchanged otherwise.

\begin{example}\label{ex:delta}
Let $(\nt,\Jt)\in\RC(\la,B)$ with $\la=\La_2$ and 
$B=(B^{1,1})^{\otimes 2}\otimes (B^{2,1})^{\otimes 3}$ of type $D_4^{(1)}$
given by
\begin{equation*}
\begin{array}[t]{r|c|c|l} \cline{2-3} 0&&&0\\ 
 \cline{2-3} 0&&\multicolumn{2}{l}{0}\\ 
 \cline{2-2}  &&\multicolumn{2}{l}{0}\\ \cline{2-2}
\end{array}
\quad
\begin{array}[t]{r|c|c|l} \cline{2-3} 0&&&0\\ \cline{2-3} &&&0\\
 \cline{2-3} 1&&\multicolumn{2}{l}{1}\\ 
 \cline{2-2} && \multicolumn{2}{l}{1}\\ \cline{2-2}
\end{array}
\quad
\begin{array}[t]{r|c|c|l} \cline{2-3} 0&&&0\\ 
 \cline{2-3} 0&&\multicolumn{2}{l}{0}\\ \cline{2-2}
\end{array}
\quad
\begin{array}[t]{r|c|c|c|l} \cline{2-4} 0&&&&0\\ \cline{2-4}
\end{array}
\end{equation*}
Here the vacancy number $p_i^{(a)}$ is written on the
left of the parts of length $i$ in $\nu^{(a)}$ and the partition
$J^{(a,i)}$ is given by the labels on the right of the parts of length
$i$ in $\nu^{(a)}$. In this case $\ell^{(1)}=\ell^{(2)}=\ell^{(3)}=1$,
$\ell^{(4)}=3$, $\lb^{(2)}=\lb^{(1)}=\infty$ and $b=\overline{3}$,
so that $\delta(\nt,\Jt)$ is
\begin{equation*}
\begin{array}[t]{r|c|c|l} \cline{2-3} 0&&&0\\ 
 \cline{2-3} 0&&\multicolumn{2}{l}{0}\\ \cline{2-2}
\end{array}
\quad
\begin{array}[t]{r|c|c|l} \cline{2-3} 0&&&0\\ \cline{2-3} &&&0\\
 \cline{2-3} 1&&\multicolumn{2}{l}{1}\\ \cline{2-2}
\end{array}
\quad
\begin{array}[t]{r|c|c|l} \cline{2-3} 1&&&0\\ \cline{2-3}
\end{array}
\quad
\begin{array}[t]{r|c|c|l} \cline{2-3} 1&&&1\\ \cline{2-3}
\end{array}.
\end{equation*}
Altogether 
\begin{equation*}
\Phi_B(\nt,\Jt)=\begin{array}{|c|}\hline \overline{3}\\ \hline \end{array}
\;\otimes\; \begin{array}{|c|}\hline \overline{4}\\ \hline \end{array}
\;\otimes\; \begin{array}{|c|}\hline 4\\ \hline 3\\ \hline \end{array}
\;\otimes\; \begin{array}{|c|}\hline \overline{1}\\ \hline 1\\ \hline \end{array}
\;\otimes\; \begin{array}{|c|}\hline 2\\ \hline 1\\ \hline \end{array}\;.
\end{equation*}
\end{example}

\subsection{Algorithm for $\ds$}\label{sec:alg spinor}
The embedding \eqref{eq:emb affine} can be extended to paths
\begin{equation}\label{eq:emb path}
\emb_{\Path}:\Path(\la,B) \hookrightarrow \Path(\lah,\Bh)
\end{equation}
where $\lah=2\la$, $B=\bigotimes_{a=1}^n (B^{a,1})^{\otimes L^a}$,
and $\Bh=\bigotimes_{a=1}^n (B^{a,1})^{\otimes \Lh^a}$ where $\Lh^a=2L^a$ 
for $1\le a\le n$. A path $b=b_L\otimes\cdots\otimes b_1\in\Path(\la,B)$
is mapped to $\emb_{\Path}(b)=\emb_B(b_L)\otimes \cdots \otimes \emb_B(b_1)$.

An analogous embedding can be defined on rigged configurations
\begin{equation}\label{eq:emb rc}
\emb_{\RC}:\RC(\la,B) \hookrightarrow \RC(\lah,\Bh)
\end{equation}
where $\emb_{\RC}(\nt,\Jt)=2(\nt,\Jt)$, meaning that all parts in $\nu^{(a)}$
and $J^{(a,i)}$ are doubled in length. 
\begin{theorem}\label{thm:emb}
For $B$ a tensor product of crystals of the form $B^{a,1}$ with $1\le a\le n-2$
the following diagram commutes:
\begin{equation*}
\begin{CD}
\RC(\la,B)        @>{\Phi_B}>> \Path(\la,B)\\
@V{\emb_{\RC}}VV @VV{\emb_B}V\\
\RC(\lah,\Bh)  @>{\Phi_{\Bh}}>> \Path(\lah,\Bh).
\end{CD}
\end{equation*}
\end{theorem}
The proof of this theorem is given in Appendix \ref{app:emb proof}.

The definition of $\ds$ is given in such a way that Theorem \ref{thm:emb}
also holds for the spinor case. Let us define
\begin{equation*}
\ds:\RC(\la,B) \to \bigcup_{\la^{\dspin}} \RC(\la^{\dspin},B')
\end{equation*}
as follows. For $(\nt,\Jt)\in \RC(\la,B)$ apply $\emb_{\RC}$ and then
a sequence of points \ref{pt:single} and \ref{pt:split}
of Theorem \ref{thm:bij} to remove the last 
tensor factor. The claim is that this rigged configuration is in the image
of $\emb_{\RC}$ so that one can apply $\emb_{\RC}^{-1}$. The result is
$\ds(\nt,\Jt)$. The proof of the well-definedness of $\ds$ is part of the proof
of Theorem~\ref{thm:bij} as given in Appendix~\ref{app:bij proof}.

\begin{example}
Let $(\nt,\Jt)\in\RC(\la,B)$ with $\la=2\La_1+\La_4$
and $B=B^{5,1}\otimes B^{2,1}\otimes (B^{1,1})^{\otimes 3}$
of type $D_5^{(1)}$ given by
\begin{equation*}
\begin{array}[t]{r|c|c|l} \cline{2-3} 2&&&1\\ \cline{2-3}
\end{array}
\quad
\begin{array}[t]{r|c|c|l} \cline{2-3} 0&&&0\\ \cline{2-3}
 0&&\multicolumn{2}{l}{0} \\ \cline{2-2}
\end{array}
\quad
\begin{array}[t]{r|c|c|l} \cline{2-3} 0&&&0\\ \cline{2-3}
 0&&\multicolumn{2}{l}{0} \\ \cline{2-2}
\end{array}
\quad
\begin{array}[t]{r|c|l} \cline{2-2} 0&&0\\ \cline{2-2}
\end{array}
\quad
\begin{array}[t]{r|c|c|l} \cline{2-3} 0&&&0\\ \cline{2-3}
\end{array}
\end{equation*}
Then $\ds(\nt,\Jt)\in \RC(\la^{\dspin},B')$ with
$\la^{\dspin}=\La_1+\La_2$ and $B'=B^{2,1}\otimes (B^{1,1})^{\otimes 3}$ is
\begin{equation*}
\begin{array}[t]{r|c|c|l} \cline{2-3} 1&&&1\\ \cline{2-3}
\end{array}
\quad
\begin{array}[t]{r|c|l} \cline{2-2} 0&&0\\ \cline{2-2}
 &&0\\ \cline{2-2}
\end{array}
\quad
\begin{array}[t]{r|c|l} \cline{2-2} 0&&0\\ \cline{2-2}
 &&0\\ \cline{2-2}
\end{array}
\quad
\begin{array}[t]{r|c|l} \cline{2-2} 0&&0\\ \cline{2-2}
\end{array}
\quad
\begin{array}[t]{r|c|l} \cline{2-2} 0&&0\\ \cline{2-2}
\end{array}
\end{equation*}
The details are given in Table \ref{tab:spinor}.
The first entry is $\emb_{\RC}(\nt,\Jt)$. The next entries are
obtained by acting with $\delta\circ \tj$. Acting with
$\emb_{\RC}^{-1}$ on the last rigged configuration yields $\ds(\nt,\Jt)$.
In the last column of the table we recorded $b$ of the algorithm $\delta$.
Hence the step in $B^{5,1}$ corresponding to this rigged configuration
is $\overline{2}\overline{5}431$.
\begin{table}
\begin{equation*}
\begin{array}{cccccc}
(\ntt,\Jtt)^{(1)}&(\ntt,\Jtt)^{(2)}&(\ntt,\Jtt)^{(3)}
&(\ntt,\Jtt)^{(4)}&(\ntt,\Jtt)^{(5)}&b\\[3mm]
\scalebox{0.8}{
\begin{array}[t]{r|c|c|c|c|l} \cline{2-5} 4&&&&&2\\ \cline{2-5}
\end{array}
}
&
\scalebox{0.8}{
\begin{array}[t]{r|c|c|c|c|l} \cline{2-5} 0&&&&&0\\ \cline{2-5}
0&&&\multicolumn{3}{l}{0}\\ \cline{2-3}
\end{array}
}
&
\scalebox{0.8}{
\begin{array}[t]{r|c|c|c|c|l} \cline{2-5} 0&&&&&0\\ \cline{2-5}
0&&&\multicolumn{3}{l}{0}\\ \cline{2-3}
\end{array}
}
&
\scalebox{0.8}{
\begin{array}[t]{l|c|c|r} \cline{2-3} 0&&&0\\ \cline{2-3}
\end{array}
}
&
\scalebox{0.8}{
\begin{array}[t]{r|c|c|c|c|l} \cline{2-5} 0&&&&&0\\ \cline{2-5}
\end{array}
}
&\\
\scalebox{0.8}{
\begin{array}[t]{r|c|c|c|c|l} \cline{2-5} 3&&&&&2\\ \cline{2-5}
\end{array}
}
&
\scalebox{0.8}{
\begin{array}[t]{r|c|c|c|l} \cline{2-4} 0&&&&0\\ \cline{2-4}
0&&&\multicolumn{2}{l}{0}\\ \cline{2-3}
\end{array}
}
&
\scalebox{0.8}{
\begin{array}[t]{r|c|c|c|l} \cline{2-4} 0&&&&0\\ \cline{2-4}
0&&&\multicolumn{2}{l}{0}\\ \cline{2-3}
\end{array}
}
&
\scalebox{0.8}{
\begin{array}[t]{l|c|c|r} \cline{2-3} 1&&&0\\ \cline{2-3}
\end{array}
}
&
\scalebox{0.8}{
\begin{array}[t]{r|c|c|c|l} \cline{2-4} 0&&&&0\\ \cline{2-4}
\end{array}
}
&
\overline{2}\\
\scalebox{0.8}{
\begin{array}[t]{r|c|c|c|c|l} \cline{2-5} 3&&&&&2\\ \cline{2-5}
\end{array}
}
&
\scalebox{0.8}{
\begin{array}[t]{r|c|c|c|l} \cline{2-4} 0&&&&0\\ \cline{2-4}
0&&&\multicolumn{2}{l}{0}\\ \cline{2-3}
\end{array}
}
&
\scalebox{0.8}{
\begin{array}[t]{r|c|c|c|l} \cline{2-4} 0&&&&0\\ \cline{2-4}
1&&&\multicolumn{2}{l}{0}\\ \cline{2-3}
\end{array}
}
&
\scalebox{0.8}{
\begin{array}[t]{l|c|c|r} \cline{2-3} 0&&&0\\ \cline{2-3}
\end{array}
}
&
\scalebox{0.8}{
\begin{array}[t]{r|c|c|l} \cline{2-3} 0&&&0\\ \cline{2-3}
\end{array}
}
&
\overline{5}\\
\scalebox{0.8}{
\begin{array}[t]{r|c|c|c|c|l} \cline{2-5} 3&&&&&2\\ \cline{2-5}
\end{array}
}
&
\scalebox{0.8}{
\begin{array}[t]{r|c|c|c|l} \cline{2-4} 0&&&&0\\ \cline{2-4}
1&&&\multicolumn{2}{l}{0}\\ \cline{2-3}
\end{array}
}
&
\scalebox{0.8}{
\begin{array}[t]{r|c|c|l} \cline{2-3} 0&&&0\\ \cline{2-3}
&&&0\\ \cline{2-3}
\end{array}
}
&
\scalebox{0.8}{
\begin{array}[t]{l|c|c|r} \cline{2-3} 0&&&0\\ \cline{2-3}
\end{array}
}
&
\scalebox{0.8}{
\begin{array}[t]{r|c|c|l} \cline{2-3} 0&&&0\\ \cline{2-3}
\end{array}
}
&
4\\
\scalebox{0.8}{
\begin{array}[t]{r|c|c|c|c|l} \cline{2-5} 3&&&&&2\\ \cline{2-5}
\end{array}
}
&
\scalebox{0.8}{
\begin{array}[t]{r|c|c|l} \cline{2-3} 0&&&0\\ \cline{2-3}
&&&0\\ \cline{2-3}
\end{array}
}
&
\scalebox{0.8}{
\begin{array}[t]{r|c|c|l} \cline{2-3} 0&&&0\\ \cline{2-3}
&&&0\\ \cline{2-3}
\end{array}
}
&
\scalebox{0.8}{
\begin{array}[t]{l|c|c|r} \cline{2-3} 0&&&0\\ \cline{2-3}
\end{array}
}
&
\scalebox{0.8}{
\begin{array}[t]{r|c|c|l} \cline{2-3} 0&&&0\\ \cline{2-3}
\end{array}
}
&
3\\
\scalebox{0.8}{
\begin{array}[t]{r|c|c|c|c|l} \cline{2-5} 2&&&&&2\\ \cline{2-5}
\end{array}
}
&
\scalebox{0.8}{
\begin{array}[t]{r|c|c|l} \cline{2-3} 0&&&0\\ \cline{2-3}
&&&0\\ \cline{2-3}
\end{array}
}
&
\scalebox{0.8}{
\begin{array}[t]{r|c|c|l} \cline{2-3} 0&&&0\\ \cline{2-3}
&&&0\\ \cline{2-3}
\end{array}
}
&
\scalebox{0.8}{
\begin{array}[t]{l|c|c|r} \cline{2-3} 0&&&0\\ \cline{2-3}
\end{array}
}
&
\scalebox{0.8}{
\begin{array}[t]{r|c|c|l} \cline{2-3} 0&&&0\\ \cline{2-3}
\end{array}
}
&
1
\end{array}
\end{equation*}
\caption{An example for $\ds$ \label{tab:spinor}}
\end{table}
\end{example}

\section{Statistics on paths and rigged configurations} \label{sec:statistics}

A statistic can be defined on both paths and rigged configurations.
In this section we define the intrinsic energy function on paths and
cocharge statistics on rigged configurations and show that the bijection $\Phi$
preserves the statistics (see Theorem \ref{thm:stat}). This gives rise to
the fermionic formula presented in Section \ref{sec:fermionic}.

In Section \ref{sec:R} we defined the combinatorial $R$-matrix
$R:B_2\otimes B_1 \to B_1\otimes B_2$.
In addition, there exists a function $H:B_2\otimes B_1\to \Z$ called the local 
energy function, that is unique up to a global additive constant.
It is constant on $J$ components and satisfies for all $b_2\in B_2$
and $b_1\in B_1$ with $R(b_2\otimes b_1)=b_1'\otimes b_2'$
\begin{equation*}
  H(e_0(b_2\otimes b_1))=
  H(b_2\otimes b_1)+
  \begin{cases}
    -1 & \text{if $\epsilon_0(b_2)>\varphi_0(b_1)$ and
    $\epsilon_0(b_1')>\varphi_0(b_2')$} \\
    1 & \text{if $\epsilon_0(b_2)\le\varphi_0(b_1)$ and
    $\epsilon_0(b_1')\le \varphi_0(b_2')$} \\
    0 & \text{otherwise.}
  \end{cases}
\end{equation*}
We shall normalize the local energy function by the condition
$H(u(B_2)\otimes u(B_1))=0$.

For a crystal $B^{k,1}$ of type $D_n^{(1)}$, the intrinsic energy
$D_{B^{k,1}}:B^{k,1}\to \Z$ is defined as follows. Let $b\in B^{k,1}$ which
is in the classical component $B(\La_{k-2j})$ (see \eqref{eq:classical comp}).
Then $D_{B^{k,1}}(b)=j$.

On the tensor product $B=B_L\otimes \cdots\otimes B_1$ of simple crystals there is an 
intrinsic energy function defined $D_B:B\to \Z$ (see for 
example~\cite[Section 2.5]{OSS:2002a})
\begin{equation} \label{eq:DNY}
  D_B = \sum_{1\le i<j\le L} H_i R_{i+1} R_{i+2}\dotsm R_{j-1}
   + \sum_{j=1}^L D_{B_j} \pi_1 R_1 R_2 \dotsm R_{j-1}.
\end{equation}
Here $R_i$ and $H_i$ denote the combinatorial $R$-matrix and local energy function
acting on the $i$-th and $(i+1)$-th tensor factors, respectively. $\pi_1$ is the
projection onto the rightmost tensor factor.

Similarly, there is a statistic on the set of rigged configurations
given by $\cc(\nt,\Jt)=\cc(\nt)+\sum_{a,i}|J^{(a,i)}|$ where
$|J^{(a,i)}|$ is the size of the partition $J^{(a,i)}$ and 
\begin{equation}\label{eq:cc}
\cc(\nt)=\frac{1}{2} \sum_{a,b\in J} \sum_{j,k\ge 1}
 (\alpha_a \mid \alpha_b) \min(j,k) m_j^{(a)}m_k^{(b)}.
\end{equation}
Let $\Phit=\Phi\circ\comp$ where $\comp:\RC(\la,B)\to\RC(\la,B)$
with $\comp(\nt,\Jt)=(\nt,\Jtt)$ is the function which complements the
riggings, meaning that $\Jtt$ is obtained from $\Jt$ by
complementing all partitions $J^{(a,i)}$ in the $m_i^{(a)}\times p_i^{(a)}$ 
rectangle.
\begin{theorem}\label{thm:stat}
Let $B=B^{k_L,1}\otimes \cdots\otimes B^{k_1,1}$ be a crystal of type
$D_n^{(1)}$ and $\la$ a dominant integral weight.
The bijection $\Phit:\RC(\la,B)\to\Path(\la,B)$ preserves the statistics,
that is $\cc(\nt,\Jt)=D(\Phit(\nt,\Jt))$ for all $(\nt,\Jt)\in
\RC(\la,B)$.
\end{theorem}
A proof of this theorem is given in Section \ref{sec:stat proof} using the results
of Sections \ref{sec:rel var maps} and \ref{sec:energy change}.

\begin{example}
For the rigged configuration of Example \ref{ex:delta} we have
\begin{equation*}
\Phit(\nt,\Jt)=\begin{array}{|c|}\hline \overline{3}\\ \hline \end{array}
\;\otimes\; \begin{array}{|c|}\hline \overline{4}\\ \hline \end{array}
\;\otimes\; \begin{array}{|c|}\hline \overline{1}\\ \hline 1\\ \hline \end{array}
\;\otimes\; \begin{array}{|c|}\hline 4\\ \hline 3\\ \hline \end{array}
\;\otimes\; \begin{array}{|c|}\hline 2\\ \hline 1\\ \hline \end{array}
\end{equation*}
and $\cc(\nt,\Jt)=D(\Phit(\nt,\Jt))=10$.
\end{example}

\subsection{Fermionic formulas} \label{sec:fermionic}

An immediate corollary of Theorems \ref{thm:bij} and \ref{thm:stat} is
the equality
\begin{equation*}
\sum_{b\in\Path(\la,B)} q^{D(b)} 
 = \sum_{(\nt,\Jt)\in \RC(\la,B)} q^{\cc(\nt,\Jt)}.
\end{equation*}
The left-hand side is in fact the one dimensional configuration sum
\begin{equation*}
X(\la,B)=\sum_{b\in\Path(\la,B)} q^{D(b)}.
\end{equation*}
The right-hand side can be simplified slightly by observing that the 
generating function of partitions in a box of width $p$ and height $m$
is the $q$-binomial coefficient
\begin{equation*}
\qbin{m+p}{m}=\frac{(q)_{p+m}}{(q)_m(q)_p},
\end{equation*}
where $(q)_m=(1-q)(1-q^2)\cdots (1-q^m)$. Hence the right-hand side becomes
the fermionic formula
\begin{equation*}
M(\la,B)=\sum_{\nt\in C(\la,B)} q^{\cc(\nt)} \prod_{\substack{i\ge 1\\ a\in J}} 
\qbin{m_i^{(a)}+p_i^{(a)}}{m_i^{(a)}}
\end{equation*}
where $m_i^{(a)}$ and $p_i^{(a)}$ are as defined in Section \ref{sec:rc}.
This proves the following result conjectured in~\cite{HKOTY:1999,HKOTT:2002}.
\begin{theorem} For $\gggg$ of type $D_n^{(1)}$, $\la$ a dominant weight and
$B=B^{k_L,1}\otimes \cdots \otimes B^{k_1,1}$
\begin{equation*}
X(\la,B)=M(\la,B).
\end{equation*}
\end{theorem}

\subsection{Relations between the various maps} \label{sec:rel var maps}

The maps $\ts,\tj,\lh,\delta$ were defined in Sections \ref{sec:var maps}
and \ref{sec:delta}. We may define similar maps using the dualities
$\dual$ and $\comp$.

For crystals let $\rh=\dual\circ\lh\circ\dual$ and
$\bs=\dual\circ\ts\circ\dual$. Explicitly, $\rh(b'\otimes b)=b'$ 
for $b'\otimes b\in B'\otimes B^{1,1}$ and
$\bs(b'\otimes m_k\cdots m_2m_1)=b'\otimes m_k\cdots m_2\otimes m_1$
for $b'\otimes m_k\cdots m_1\in B'\otimes \Bh^{k,1}$ for $2\le k\le n-1$.
Similarly $\bs(b'\otimes m_n\cdots m_2m_1)=b'\otimes m_n\cdots m_2\otimes m_1$
for $b'\otimes m_n\cdots m_1\in B'\otimes E^{n-1,1}$ or $E^{n,1}$.

On the rigged configuration side define $\dt=\comp\circ\delta\circ\comp$
and $\bj=\comp\circ\tj\circ\comp$. Explicitly, $\bj$ adds parts of lengths one to the
same rigged partitions as $\tj$; however in this case they are not
singular, but of label 0.

In the following $[\cdot,\cdot]$ denotes the commutator.

\begin{lemma} \label{lem:split-delta}
We have $[\bj,\delta]=0$ and $[\tj,\dt]=0$.
\end{lemma}
\begin{proof}
The relation $[\tj,\dt]=0$ follows from $[\bj,\delta]=0$ by conjugation by $\comp$.
Hence it suffices to prove $[\bj,\delta]=0$.

Let $(\nt,\Jt)\in\RC(\la,B)$ where $B=\Bh^{1,1}\otimes \Bh^{k,1}\otimes B'$. We assume throughout
the proof that $1\le k<n$. The cases $B=\Bh^{1,1}\otimes E^{n,1}\otimes B'$ and 
$B=\Bh^{1,1}\otimes E^{n-1,1}\otimes B'$ can be proven in a similar fashion. 
Let $\ell^{(a)}$ and $\lb^{(a)}$ (resp. $s^{(a)}$ and $\sbar^{(a)}$) be the singular strings 
selected by $\delta$ in $(\nt,\Jt)$ (resp. $\bj(\nt,\Jt)$). Clearly, $s^{(a)}\le \ell^{(a)}$ 
and $\sbar^{(a)}\le \lb^{(a)}$ since $\bj$ changes neither vacancy numbers nor riggings by 
Lemma \ref{lem:tj} and its definition $\bj=\comp\circ \tj\circ\comp$. If $s^{(a)}=\ell^{(a)}$ 
and $\sbar^{(a)}=\lb^{(a)}$ for all $a$, the lemma holds.

Assume that $b$ is minimal such that $s^{(b)}<\ell^{(b)}$. Certainly we must have
$1\le b<k$, $s^{(a)}=\ell^{(a)}=1$ for $1\le a<b$, $s^{(b)}=1<\ell^{(b)}$,
$m_1^{(b)}=0$ and $p_1^{(b)}=0$. One may easily verify that
(see also \cite[Eq. (3.10)]{KSS:2002})
\begin{equation*}
\begin{split}
-p_{i-1}^{(a)}+2p_i^{(a)}-p_{i+1}^{(a)}
=&-\sum_{b\in J} \left( \alpha_a \mid \alpha_b\right) m_i^{(b)}+\delta_{i,1} L^a\\
\ge & -\sum_{b\in J} \left( \alpha_a \mid \alpha_b\right) m_i^{(b)},
\end{split}
\end{equation*}
where $\delta_{i,j}=1$ if $i=j$ and $0$ else. Using $m_1^{(b)}=0$ and $p_1^{(b)}=0$
this reads for $i=1$ and $a=b$
\begin{equation*}
-p_2^{(b)}=m_1^{(b-1)}+m_1^{(b+1)}+L^b.
\end{equation*}
Since $p_i^{(a)}\ge 0$ for all admissible rigged configurations, this implies that
all quantities in this equation must be zero. For $b>1$ the condition $m_1^{(b-1)}=0$
contradicts $\ell^{(b-1)}=1$. For $b=1$ we have $L^1\ge 1$ which makes $p_2^{(1)}\ge 0$ impossible.
Hence $s^{(a)}=\ell^{(a)}$ for all $a$ as desired. The proof that $\sbar^{(a)}=\lb^{(a)}$ is
analogous.
\end{proof}

\begin{lemma} \label{lem:delta-delta}
We have $[\delta,\dt]=0$.
\end{lemma}
The proof of Lemma \ref{lem:delta-delta} is outlined in Appendix \ref{app:delta-delta}.

We say that $f$ corresponds to $g$ under $\Phi$ if the following
diagram commutes:
\begin{equation*}
\begin{CD}
 \RC(\la,B)  @>{\Phi_B}>> \Path(\la,B)\\ 
 @V{f}VV @VV{g}V\\
 \RC(\la',B') @>>{\Phi_{B'}}> \Path(\la',B').
\end{CD}
\end{equation*}

\begin{prop}\label{prop:corresp}
The following maps correspond under $\Phi$:
\begin{enumerate}
\item $\bj$ and $\bs$
\item $\tj$ and $\ts$
\item $\delta$ and $\lh$
\item $\dt$ and $\rh$
\item $\comp$ and $\dual$
\item $\id$ and $R$, where $R=R_{i_1}\circ\cdots\circ R_{i_k}$ is a composition of
the $R_i$'s.
\end{enumerate}
\end{prop}
\begin{proof}
Points (2) and (3) hold by definition. 

We prove (1) by induction on $B$. For simplicity we drop the dependency on $\lambda$
in all diagrams to follow. Also for brevity we denote 
$\Bh^{k,1}\otimes B \otimes \Bh^{k',1}$ symbolically by $\vdom B\vdom$.
The splitting $\Bh^{k,1}\to \Bh^{k-1,1}\otimes \Bh^{1,1}$ will be denoted
by $\vdom \to \sq\;\sq$ and so on. Consider the following diagram:
\begin{equation}
\label{eq:dia1-1}
\xymatrix{
{\RC(\vdom B\vdom)} \ar[rrr]^{\Phi} \ar[ddd]_{\bj} \ar[dr]^{\tj}
& & & {\Path(\vdom B\vdom)} \ar[ddd]^{\bs} \ar[dl]_{\ts} \\
 & {\RC(\sq\;\sq B\vdom)} \ar[r]^{\Phi} \ar[d]_{\bj} &
	 {\Path(\sq\;\sq B\vdom)} \ar[d]^{\bs} & \\
 & {\RC(\sq\;\sq B\sq\;\sq)} \ar[r]_{\Phi} & {\Path(\sq\;\sq B\sq\;\sq)} & \\
{\RC(\vdom B\sq\;\sq)} \ar[ur]_{\tj} \ar[rrr]_{\Phi} & & & {\Path(\vdom B\sq\;\sq)}
	\ar[ul]^{\ts}
}
\end{equation}
This diagram can be viewed as a cube with the front face given by the large square.
The back face commutes by induction. It is easy to check that the left and right
faces commute. The top and bottom face commute by (2).
Hence it follows from \cite[Lemma 5.3]{KSS:2002} that the front face commutes 
since $\ts$ is injective.

Next consider the diagram:
\begin{equation}
\label{eq:dia1-2}
\xymatrix{
{\RC(\sq B\vdom)} \ar[rrr]^{\Phi} \ar[ddd]_{\bj} \ar[dr]^{\delta}
& & & {\Path(\sq B\vdom)} \ar[ddd]^{\bs} \ar[dl]_{\lh} \\
 & {\RC(B\vdom)} \ar[r]^{\Phi} \ar[d]_{\bj} &
	 {\Path(B\vdom)} \ar[d]^{\bs} & \\
 & {\RC(B\sq\;\sq)} \ar[r]_{\Phi} & {\Path(B\sq\;\sq)} & \\
{\RC(\sq B\sq\;\sq)} \ar[ur]_{\delta} \ar[rrr]_{\Phi} & & & {\Path(\sq B\sq\;\sq)}
	\ar[ul]^{\lh}
}
\end{equation}
The back face commutes by induction, the top and bottom faces commute by the
definition of $\Phi$. It is easy to check that the right face commutes.
The left face commutes by Lemma \ref{lem:split-delta}. Hence the front
face commutes by \cite[Lemma 5.3]{KSS:2002}. Diagrams \eqref{eq:dia1-1}
and \eqref{eq:dia1-2} together prove (1).

For the proof of (4) we use the same abreviated notation as before.
$\Bh^{1,1}$ is symbolically denoted by $\sq$. Consider the following diagram:
\begin{equation}
\label{eq:dia4-1}
\xymatrix{
{\RC(\sq B\sq)} \ar[rrr]^{\Phi} \ar[ddd]_{\dt} \ar[dr]^{\delta}
& & & {\Path(\sq B\sq)} \ar[ddd]^{\rh} \ar[dl]_{\lh} \\
 & {\RC(B\sq)} \ar[r]^{\Phi} \ar[d]_{\dt} &
	 {\Path(B\sq)} \ar[d]^{\rh} & \\
 & {\RC(B)} \ar[r]_{\Phi} & {\Path(B)} & \\
{\RC(\sq B)} \ar[ur]_{\delta} \ar[rrr]_{\Phi} & & & {\Path(\sq B)}
	\ar[ul]^{\lh}
}
\end{equation}
The back face commutes by induction and the top and bottom face commute by (3).
It is easy to verify that the right face commutes. The left face commutes by
Lemma \ref{lem:delta-delta}. Hence the front face commutes up to $\delta$.
Since $\delta$ is injective, it suffices to show that both ways around the front
face result in elements with the same highest weight. This follows 
from~\cite[Propositions 5.11 and 8.5]{SS}.

Next consider the diagram:
\begin{equation}
\label{eq:dia4-2}
\xymatrix{
{\RC(\vdom B\sq)} \ar[rrr]^{\Phi} \ar[ddd]_{\dt} \ar[dr]^{\tj}
& & & {\Path(\vdom B\sq)} \ar[ddd]^{\rh} \ar[dl]_{\ts} \\
 & {\RC(\sq\;\sq B\sq)} \ar[r]^{\Phi} \ar[d]_{\dt} &
	 {\Path(\sq\;\sq B\sq)} \ar[d]^{\rh} & \\
 & {\RC(\sq\;\sq B)} \ar[r]_{\Phi} & {\Path(\sq\;\sq B)} & \\
{\RC(\vdom B)} \ar[ur]_{\tj} \ar[rrr]_{\Phi} & & & {\Path(\vdom B)}
	\ar[ul]^{\ts}
}
\end{equation}
The back face commutes by induction, the top and bottom face commute by (2).
It is easy to verify that the right face commutes. The left face commutes by
Lemma \ref{lem:split-delta}. Hence the front face commutes, since $\ts$ is injective.
Diagrams \eqref{eq:dia4-1} and \eqref{eq:dia4-2} together prove that (4) holds.

For the proof of point (5) consider the following diagram:
\begin{equation}
\label{eq:dia5-1}
\xymatrix{
{\RC(\sq B)} \ar[rrr]^{\Phi} \ar[ddd]_{\comp} \ar[dr]^{\delta}
& & & {\Path(\sq B)} \ar[ddd]^{\dual} \ar[dl]_{\lh} \\
 & {\RC(B)} \ar[r]^{\Phi} \ar[d]_{\comp} &
	 {\Path(B)} \ar[d]^{\dual} & \\
 & {\RC(B^\dual)} \ar[r]_{\Phi} & {\Path(B^\dual)} & \\
{\RC(B^\dual\sq)} \ar[ur]_{\dt} \ar[rrr]_{\Phi} & & & {\Path(B^\dual\sq)}
	\ar[ul]^{\rh}
}
\end{equation}
The back face commutes by induction. The left and right faces commute by
definition. The top and bottom faces commute by (3) and (4), respectively.
Since $\dt$ is injective, the front face commutes by \cite[Lemma 5.3]{KSS:2002}.

Next consider the following diagram:
\begin{equation}
\label{eq:dia5-2}
\xymatrix{
{\RC(\vdom B)} \ar[rrr]^{\Phi} \ar[ddd]_{\comp} \ar[dr]^{\tj}
& & & {\Path(\vdom B)} \ar[ddd]^{\dual} \ar[dl]_{\ts} \\
 & {\RC(\sq\;\sq B)} \ar[r]^{\Phi} \ar[d]_{\comp} &
	 {\Path(\sq\;\sq B)} \ar[d]^{\dual} & \\
 & {\RC(B^\dual\sq\;\sq)} \ar[r]_{\Phi} & {\Path(B^\dual \sq\;\sq)} & \\
{\RC(B^\dual\vdom)} \ar[ur]_{\bj} \ar[rrr]_{\Phi} & & & {\Path(B^\dual\vdom)}
	\ar[ul]^{\bs}
}
\end{equation}
The back face commutes by induction. The left and right faces commute by
definition. The bottom and top faces commute by (1) and (2), respectively.
Since $\bj$ is injective, the front face commutes by \cite[Lemma 5.3]{KSS:2002}.
Hence (5) holds by diagrams \eqref{eq:dia5-1} and \eqref{eq:dia5-2}.

For the proof of (6) we may assume that $R=R_i$ is the combinatorial $R$-matrix
acting on the $i$-th and $(i+1)$-th tensor factor. By induction we may assume that
$i=L-1$ where $L$ is the number of tensor factors in $B$. By part (5) we may
assume that $i=1$. Again by induction we may assume that $L=2$, which means that
$B=B^{k_2,1}\otimes B^{k_1,1}$.
Hence, since $R$ is a crystal isomorphism, it suffices to check (6) for highest weight
elements in $B=B^{k',1}\otimes B^{k,1}$. The combinatorial $R$-matrix
$R:B^{k',1}\otimes B^{k,1}\to B^{k,1}\otimes B^{k',1}$ for all highest weight elements
is given explicitly in Appendix \ref{app:R matrix}. It is tedious but straightforward 
to check (6) for these.
\end{proof}

\subsection{Properties of the statistics} \label{sec:energy change}

In this section we show how the statistics $\cc$ and $D$ change under 
$\tj$ and $\bs$, respectively. These results will be used in the next
section in the proof of Theorem \ref{thm:stat}.
For the next two lemmas we assume that $B$ is the tensor product of 
crystals $\Bh^{k,1}$ for $1\le k<n$, $E^{n-1,1}$ or $E^{n,1}$.
\begin{lemma}\label{lem:cc}
For $(\nt,\Jt)\in\RC(\la,B)$ we have
\begin{equation*}
\cc(\tj(\nt,\Jt))-\cc(\nt,\Jt)=
\begin{cases}
 (\epsilon_1-\epsilon_k \mid \vL)+1 
 & \text{if $B=\Bh^{k,1}\otimes B'$ for $1\le k<n$}\\
 (\epsilon_1-\epsilon_n \mid \vL)+1 & \text{if $B=E^{n,1}\otimes B'$}\\
 (\epsilon_1+\epsilon_n \mid \vL)+1 & \text{if $B=E^{n-1,1}\otimes B'$.}
\end{cases}
\end{equation*}
\end{lemma}
\begin{proof}
First consider the case $B=\Bh^{k,1}\otimes B'$. Then by \eqref{eq:cc}
\begin{equation*}
\begin{aligned}
\cc(\tj(\nt))=&\frac{1}{2}\sum_{a,b}\sum_{i,j} 
 (\alpha_a\mid \alpha_b)\min(i,j)\\
 &\times \Bigl(m_i^{(a)}+\chi(i=1)\chi(a<k)\Bigr)
         \Bigl(m_j^{(b)}+\chi(j=1)\chi(b<k)\Bigr)\\
=&\cc(\nt)+\sum_{b,j}(\alpha_1+\cdots+\alpha_{k-1}\mid \alpha_b) m_j^{(b)}+1.
\end{aligned}
\end{equation*}
Since $\tj$ adds a singular string of length one to the first $k-1$ rigged
partitions, the change from the labels is, using \eqref{eq:vac},
\begin{equation*}
\begin{split}
p_1^{(1)}+p_1^{(2)}+\cdots+p_1^{(k-1)}
=(\alpha_1+\cdots+\alpha_{k-1}\mid \vL-\sum_{b,j}\alpha_b m_j^{(b)}).
\end{split}
\end{equation*}
Hence
\begin{equation*}
\cc(\tj(\nt,\Jt))-\cc(\nt,\Jt)=(\epsilon_1-\epsilon_k\mid \vL)+1.
\end{equation*}
The cases $B=E^{n,1}\otimes B'$ and $B=E^{n-1,1}\otimes B'$ are similar.
\end{proof}

\begin{lemma}\label{lem:D}
For $p\in\Path(\la,B)$ we have
\begin{equation*}
D(\bs(p))-D(p)=
\begin{cases}
 (\epsilon_1-\epsilon_k \mid \vL)+1 
 & \text{if $B=B'\otimes \Bh^{k,1}$ for $1\le k<n$}\\
 (\epsilon_1-\epsilon_n \mid \vL)+1 & \text{if $B=B'\otimes E^{n,1}$}\\
 (\epsilon_1+\epsilon_n \mid \vL)+1 & \text{if $B=B'\otimes E^{n-1,1}$.}
\end{cases}
\end{equation*}
\end{lemma}
\begin{proof}
It follows from the structure of \eqref{eq:DNY} that we can restrict ourselves
to the case when $B$ consists of one or two tensor factors. Let $b\in \Bh^{k,1}$ 
($E^{n-1,1}$ or $E^{n,1}$) and $b'\in \Bh^{k',1}$ ($E^{n-1,1}$ or $E^{n,1}$) such 
that $e_i(b)=0$ and $e_i(b'\otimes b)=0$ for $1\le i\le n$.
Furthermore, write $\bs(b)=b_2\otimes b_1\in \Bh^{k-1,1}\otimes B^{1,1}$ 
($\Bh^{n-1,1}\otimes B^{1,1}$) and
\begin{equation*}
\begin{split}
R(b'\otimes b) &= \bh \otimes \bh'\\
R(b'\otimes b_2) &= \bullet \otimes \bt'\\
R(\bt'\otimes b_1) &= \bullet \otimes \btt', 
\end{split}
\end{equation*}
where the crystal elements denoted $\bullet$ are not of interest to us here.
Let $\vL'$ be the weight corresponding to $\Bh^{k',1}$ (resp. $E^{n-1,1}$
or $E^{n,1}$) (see \eqref{eq:L}). Then it suffices to prove the following two 
equations:
\begin{align}
\label{eq:1tensor}
&D(b_2\otimes b_1)-D(b) = 1\\
\label{eq:2tensor}
&H(b'\otimes b_2)+H(\bt'\otimes b_1)+D(\btt')-H(b'\otimes b)-D(\bh')\\
&\notag\qquad\qquad\qquad=
\begin{cases} (\epsilon_1-\epsilon_k\mid \vL') & \text{if $b\in \Bh^{k,1}$}\\
              (\epsilon_1-\epsilon_n\mid \vL') & \text{if $b\in E^{n,1}$}\\
              (\epsilon_1+\epsilon_n\mid \vL') & \text{if $b\in E^{n-1,1}$}.
\end{cases}
\end{align}

To prove \eqref{eq:1tensor} and \eqref{eq:2tensor} we use the same notation as
in Appendix \ref{app:R matrix}. Namely, a column with entries
$\cdots a_4\cdots a_3 a_2 \cdots a_1$, where $a_2\cdots a_1$ denotes a
sequence of consecutive integers $a_2\, a_2-1\, a_2-2\,\ldots a_1+1\, a_1$, is compactly 
written as $\left[a_1..a_2 \mid a_3..a_4 \mid \cdots \right]$.

Let us first concentrate on \eqref{eq:1tensor} for $b\in \Bh^{k,1}$. 
Since $e_i(b)=0$ for $1\le i\le n$, $b$ is of the form 
$b=[1..p \mid \ol{p}..\ol{\ell}]$. Then $b_2\otimes b_1=[2..p \mid \ol{p}..\ol{\ell}] 
\otimes [1]$ and $R(b_2\otimes b_1)=\bt_1\otimes \bt_2=[\ell-1]\otimes [1..p-1 \mid 
\ol{p-1}..\ol{\ell-1}]$. Hence $D(b)=p-\ell+1$ and $D(b_2\otimes b_1)
=H(b_2\otimes b_1)+D(b_1)+D(\bt_2)=1+0+(p-\ell+1)=p-\ell+2$ which 
proves \eqref{eq:1tensor}. The highest weight vector in $E^{n,1}$ is
$b=[1..n]$ and $b_2\otimes b_1=[2..n]\otimes [1]$, so
that $R(b_2\otimes b_1)=\bt_1\otimes \bt_2=[n]\otimes [1..n-1]$.
In this case $D(b)=0$ and $D(b_2\otimes b_1)=H(b_2\otimes b_1)+D(b_1)+D(\bt_2)
=1+0+0=1$. The highest weight vector in $E^{n-1,1}$ is
$b=[1..n-1\, \ol{n}]$ and the calculation is analogous to the case $E^{n,1}$.

Equation \eqref{eq:2tensor} can be proven for each case listed in
Appendix \ref{app:R matrix}, where it should be noted that the results for
$B^{k',1}\otimes B^{k,1}$ with $1\le k'\le k\le n-2$ can be extended to
$\Bh^{k',1}\otimes \Bh^{k,1}$ with $1\le k'\le k\le n-1$ and
$\Bh^{k',1}\otimes \Bh^{n,1}$ ($\Bh^{k',1}\otimes \Bbh^{n,1}$) with $k'<n$. 

Here we show the details for Case 1 for
$k'<k\le n$ and $r\le \ell$. In this case 
$b'\otimes b = [1..\ell \mid k+1..p \mid \ol{p}..\ol{q} \mid \ol{\ell}..\ol{r}]
\otimes [1..k]$, $b_2=[2..k]$ and $b_1=[1]$. Also by the results of 
Appendix \ref{app:R matrix}, $H(b'\otimes b)=k'-\ell$ and $D(\bh')=0$.
To determine $H(b'\otimes b_2)$, $H(\bt'\otimes b_1)$ and $D(\btt')$,
note that
\begin{equation*}
e_q e_{q+1} \cdots e_{n-2}\; e_n e_{n-1} \cdots e_1 (b'\otimes b_2)
=[1..\ell \mid k..p \mid \ol{p}..\ol{q+1} \mid \ol{\ell}..\ol{r}]
\otimes [1..k-1].
\end{equation*}
Acting on this by the $R$-matrix one obtains 
\begin{equation*}
\begin{cases}
[1..\ell+1 \mid k'+1..k-2 \mid \ol{k-2}..\ol{q+1} \mid \ol{\ell+1}..\ol{r}]
\otimes [1..k'] & \text{if $p=k-1$, $k-q$ even}\\
[1..\ell \mid k'+1..p \mid \ol{p}..\ol{q+1} \mid \ol{\ell}..\ol{r}]
\otimes [1..k'] & \text{else}
\end{cases}
\end{equation*}
so that in both cases $\bt'=[2..k'+1]$ and $\btt'=[1..k']$.
Since $H$ is constant on classical components, $H(b'\otimes b_2)=k'-\ell$,
$H(\bt'\otimes b_1)=1$ and $D(\btt')=0$. Hence the result of \eqref{eq:2tensor} is
1 as claimed for $k'<k\le n$.

In Appendix \ref{app:R matrix}, the $R$-matrix is only worked out for the 
case $k'\le k$. To prove \eqref{eq:2tensor} for $k\le k'$, consider the equation
under the duality operation $\dual$ which swaps tensor factors.
Since $\ts=\dual\circ\bs\circ\dual$, \eqref{eq:2tensor} becomes
\begin{equation}\label{eq:2tensor*}
\begin{split}
H(b_2'\otimes b)&+H(b_1'\otimes \bt)+D(\btt)-H(b'\otimes b)-D(\bh)\\
&=
\begin{cases} (\epsilon_1-\epsilon_{k'}\mid \vL) & \text{if $b'\in \Bh^{k',1}$}\\
              (\epsilon_1-\epsilon_n\mid \vL) & \text{if $b'\in E^{n,1}$}\\
              (\epsilon_1+\epsilon_n\mid \vL) & \text{if $b'\in E^{n-1,1}$}
\end{cases}
\end{split}
\end{equation}
where $\ts(b')=b_1'\otimes b_2'$, $R(b_2'\otimes b)=\bt\otimes \bullet$,
$R(b_1'\otimes \bt)=\btt\otimes \bullet$, and $\vL$ is the weight corresponding to 
$\Bh^{k,1}$ (resp. $E^{n,1}$, $E^{n-1,1}$). Again \eqref{eq:2tensor*} can be
verified for all cases of Appendix \ref{app:R matrix}. Here we demonstrate
that it holds for Case 1 with $k'\le k<n$ and $r\le \ell$, so that
$b'\otimes b = [1..\ell \mid k+1..p \mid \ol{p}..\ol{q} \mid \ol{\ell}..\ol{r}]
\otimes [1..k]$, $b_1'=[\ol{r}]$ and $b_2'=[1..\ell \mid k+1..p \mid \ol{p}..\ol{q}
\mid \ol{\ell}..\ol{r+1}]$. Hence $H(b'\otimes b)=k'-\ell$ and 
$H(b_2'\otimes b)=k'-\ell-1$. Next let us calculate $D(\bh)$.
By Appendix \ref{app:R matrix}, $\bh=[1..\ell \mid k'+1..p \mid \ol{p}..\ol{q}
\mid \ol{\ell}..\ol{r}]$ if $p\neq k$ or $k-q$ odd. Acting with a product
of $e_i$ with $i\neq 0$ to make $\bh$ into a highest weight vector,
one finds that $D(\bh)=d+\max\{0,p-q+1-k'+\ell\}$. However, since
$k'=2\ell+2p-k-q-r+2$ we have $p-q+1-k'+\ell=(k-p)+(r-\ell)-1<0$
because $k-p\le 0$ and $r-\ell\le 0$. This implies that
$D(\bh)=d$. Similarly, for $p=k$ and $k-q$ even one finds that
$D(\bh)=d+1$. By Appendix \ref{app:R matrix}
\begin{equation*}
\bt= \begin{cases}
[1..\ell+1 \mid k'..k-1 \mid \ol{k-1}..\ol{q} \mid \ol{\ell+1}..\ol{r+1}]
& \text{if $p=k$, $k-q$ even}\\
[1..\ell \mid k'..p \mid \ol{p}..\ol{q} \mid \ol{\ell}..\ol{r+1}]
& \text{else.}
\end{cases}
\end{equation*}
In both cases the highest weight vector obtained from $b_1'\otimes \bt$
by the application of a product of $e_i$'s with $i\neq 0$ is of the form
$[\ol{y}] \otimes [1..x \mid \ol{x}..\ol{y+1}]$ where $x=k-d$ when $p=k$ and 
$k-q$ even and $x=k-d+1$ otherwise. The variable $y$ is fixed by the requirement
that $k=2x-y$. By Case 5a of Appendix \ref{app:R matrix} we have
\begin{equation*}
R([\ol{y}]\otimes [1..x\mid \ol{x}..\ol{y+1}]) = 
[1..x-1 \mid \ol{x-1}..\ol{y}\mid \ol{1}] \otimes [1]
\end{equation*}
so that $D(\btt)=x-y+1=k+1-x$ which is $d+1$ if $p=k$ and $k-q$ even and 
$d$ otherwise. Hence altogether the left-hand side of \eqref{eq:2tensor*} yields
0 as expected.

Finally let us show that \eqref{eq:2tensor*} yields 2 for the tensor product
$E^{n-1,1}\otimes E^{n,1}$ and 0 for $E^{n,1}\otimes E^{n,1}$.
The highest weight elements in $E^{n-1,1}\otimes E^{n,1}$ are of the form
$b'\otimes b=[1..n-2\ell-1 \mid \ol{n}..\ol{n-2\ell}]\otimes [1..n]$. 
By Appendix \ref{app:R matrix}, $H(b'\otimes b)=2\ell$ and $D(\bh)=0$.
By the definition of $\ts$, $b_1'=[\ol{n-2\ell}]$ and 
$b_2'=[1..n-2\ell-1 \mid \ol{n}..\ol{n-2\ell+1}]$. Furthermore
\begin{equation*}
\bt\otimes \bullet = R(b_2'\otimes b)=[1..n-2\ell-1 \mid n \ol{n}..\ol{n-2\ell+1}]
\otimes [1..n-1]
\end{equation*}
with $H(b_2'\otimes b)=2\ell$. The highest weight vector corresponding to
$b_1'\otimes \bt$ is $[\ol{n}]\otimes [1..n]$ with
$R([\ol{n}]\otimes [1..n])=[1..n-1\mid \ol{1}]\otimes [1]$. Hence $H(b_1'\otimes \bt)=1$
and $D(\btt)=1$, so that \eqref{eq:2tensor*} yields indeed 2.

The highest weight elements in $E^{n,1}\otimes E^{n,1}$ are of the form
$b'\otimes b=[1..n-2\ell \mid \ol{n}..\ol{n-2\ell+1}]\otimes [1..n]$. 
By Appendix \ref{app:R matrix}, $H(b'\otimes b)=2\ell$ and $D(\bh)=0$.
By the definition of $\ts$, $b_1'=[\ol{n-2\ell+1}]$ and 
$b_2'=[1..n-2\ell \mid \ol{n}..\ol{n-2\ell+2}]$. Furthermore
\begin{equation*}
\bt\otimes \bullet = R(b_2'\otimes b)=[1..n-2\ell+1 \mid \ol{n-1}..\ol{n-2\ell+1}]
\otimes [1..n-1]
\end{equation*}
with $H(b_2'\otimes b)=2\ell-1$. The highest weight vector corresponding to
$b_1'\otimes \bt$ is $[1]\otimes [1..n-1 \mid \ol{n-1}]$ with
$R([1]\otimes [1..n-1 \mid \ol{n-1}])=[1..n-1\mid \ol{n-1}]\otimes [1]$. Hence 
$H(b_1'\otimes \bt)=0$ and $D(\btt)=1$, so that \eqref{eq:2tensor*} yields indeed 0.

All other cases can be shown in a similar fashion.
\end{proof}

\subsection{Proof of Theorem \ref{thm:stat}} \label{sec:stat proof}

By Proposition \ref{prop:corresp} the following diagram commutes
\begin{equation*}
\begin{CD}
\RC(\la,B) @>{\Phi}>> \Path(\la,B) @>{\dual}>> \Path(\la,B)\\
@V{\tj}VV @VV{\ts}V @VV{\bs}V\\
\RC(\la,\Bt) @>>{\Phi}> \Path(\la,\Bt) @>>{\dual}> \Path(\la,\Bt).
\end{CD}
\end{equation*}
Since $\Phit=\dual\circ\Phi$, this implies that $\tj$ corresponds to
$\bs$ under $\Phit$.

Let $(\nt,\Jt)\in\RC(\la,B)$ and $b=\Phit(\nt,\Jt)\in\Path(\la,B)$.
Using \eqref{eq:emb path} and \eqref{eq:emb rc} define
$(\nht,\Jht)=\emb_{\RC}(\nt,\Jt)\in\RC(\lah,\Bh)$ and 
$\bh=\emb_{\Path}(b)\in\Path(\lah,\Bh)$ with $\lah$ and $\Bh$ as defined in 
Section \ref{sec:alg spinor}.
It is easy to see that $\cc(\nht,\Jht)=2\cc(\nt,\Jt)$ since all parts of
the partitions and all riggings double. Similarly, by \cite[Eq. (6.23)]{OSS}
it follows that $D_{\Bh}(\bh)=2D_B(b)$. Since the embedding map commutes with
the bijection $\Phit$ by Theorem \ref{thm:emb} and the definition of the spinor
case, we may restrict ourselves to the image of $\emb$.

There exists a sequence $S_{\Path}$ of maps $\bs$ and $R$
which transforms $\bh$ into a path of single boxes, that is 
$S_{\Path}:\Path(\lah,\Bh)\to\Path(\lah,(B^{1,1})^N)$ where 
$N=\sum_{a=1}^n a \Lh^a + n\Lbh^n$.
By Proposition \ref{prop:corresp} there is a corresponding sequence 
$S_{\RC}:\RC(\lah,\Bh)\to \RC(\lah,(B^{1,1})^N)$ under $\Phit$ consisting of maps
$\tj$ and $\id$.

By Lemmas \ref{lem:cc} and \ref{lem:D} and the fact that $R$ and $\id$ do not
change the statistics, it follows that for $\bh=\Phit(\nht,\Jht)$
\begin{equation*}
D(S_{\Path}(\bh)) = \cc(S_{\RC}(\nht,\Jht)) \qquad \text{implies that} \qquad
D(\bh) = \cc(\nht,\Jht),
\end{equation*}
which in turn implies that $D(b)=\cc(\nt,\Jt)$.
But Theorem \ref{thm:stat} for the case $B=(B^{1,1})^N$ has already been
proven in \cite{OSS:2002}.

\appendix

\section{Proof of Theorem \ref{thm:bij}}
\label{app:bij proof}

The proof proceeds by induction on $B$. Given a tensor product
$B$ of tensor factors of the form $\Bh^{k,1}$ for $1\le k\le n-1$ or
$B^{k,1}$ for $k=n-1,n$, there is a unique sequence of transformations 
$-:B=B^{1,1}\otimes B'\to B'$, $\ts:B=\Bh^{k,1}\otimes B'\to 
B^{1,1}\otimes B^{k-1,1}\otimes B'$ for $2\le k\le n-1$ and
$-_s:B=B^{k,1}\otimes B'\to B'$ for $k=n-1,n$ resulting in the empty
crystal. If $B$ is empty then both $\RC(\la,B)$ and $\Path(\la,B)$
are the empty set unless $\la=\es$, in which case $\Path(\la,B)$
is the singleton consisting of the empty path, $\RC(\la,B)$ is the 
singleton consisting of the empty rigged configuration, and 
$\Phi_\es$ is the unique bijection $\RC(\es,\es)
\to \Path(\es,\es)$.

\subsection*{Proof of \eqref{pt:single}}
Any map $\Phi_B$ satisfying \eqref{pt:single} is injective by
definition and unique by induction. It was proven in \cite{OSS:2002}
that $-$ and $\delta$ are bijections.

\subsection*{Proof of \eqref{pt:split}}
To prove the existence and uniqueness of a bijection satisfying
\eqref{pt:split}, it must be shown that $\Phi_{\Bt}$ restricts to
a bijection $\Phi_{\Bt}:\tj(\RC(\la,B))\to \ts(\Path(\la,B))$.
The proof relies on two lemmas for which we need the following notation.
Let $(\nt,\Jt)\in \tj(\RC(\la,B))$ and $(\ntt,\Jtt)=\tj\circ\delta(\nt,\Jt)$.
Let $\ell^{(a)},\lb^{(a)}$ be the length of the selected singular
strings in $(\nt,\Jt)$ by $\delta$, and $\lt^{(a)},\lbt^{(a)}$ 
the ones selected in $(\ntt,\Jtt)$ by $\delta$.

\begin{lemma} \label{lem:ell}
The following inequalities hold:
\begin{enumerate}
\item[(1)] $\lt^{(a)} \ge \ell^{(a+1)}$ for $1\le a \le n-3$
\item[(2)] $\lt^{(n-2)} \ge \min(\ell^{(n-1)},\ell^{(n)})$
\item[(3)] $\lt^{(n-1)} \ge \lb^{(n-2)}$ 
           if $\ell^{(n-1)}=\min(\ell^{(n-1)},\ell^{(n)})$\\
           $\lt^{(n)} \ge \lb^{(n-2)}$ if $\ell^{(n)}=\min(\ell^{(n-1)},\ell^{(n)})$
\item[(4)] $\lbt^{(a)} \ge \lb^{(a-1)}$ for $2\le a\le n-1$.
\end{enumerate}
\end{lemma}

Similarly, let $b\in \Path(\la,B)$, $\bt=-\circ\ts(b)$ and $\btt=-\circ\ts(\bt)$. 
More precisely, if $b=m_k m_{k-1} \cdots m_1 \otimes b'$, 
then $\bt=m_{k-1} m_{k-2}\cdots m_1 \otimes b'$ and 
$\btt=m_{k-2}\cdots m_1 \otimes b'$.
By induction $(\ntt,\Jtt)=\Phi^{-1}(\bt)$ and $(\nttt,\Jttt)=\Phi^{-1}(\btt)$ exist. 
Let $\st^{(a)},\sbt^{(a)}$ be the lengths of the selected strings by $\delta^{-1}$ 
in $(\nttt,\Jttt)$ adding $m_{k-1}\in B^{1,1}$ to the path, and $
s^{(a)},\sbar^{(a)}$ the length of the 
selected strings in $(\ntt,\Jtt)$ adding the step $m_k\in B^{1,1}$ to the path.

\begin{lemma} \label{lem:s}
The following inequalities hold:
\begin{enumerate}
\item[(1)] $\sbar^{(a)}\le \sbt^{(a+1)}$ for $1\le a \le n-3$
\item[(2)] $\sbar^{(n-2)} \le \max(\st^{(n-1)},\st^{(n)})$
\item[(3)] $s^{(n-1)}\le \st^{(n-2)}$ if $\st^{(n-1)}=\max(\st^{(n-1)},\st^{(n)})$\\
           $s^{(n)}\le \st^{(n-2)}$ if $\st^{(n)}=\max(\st^{(n-1)},\st^{(n)})$
\item[(4)] $s^{(a+1)}\le \st^{(a)}$ for $0\le a \le n-3$.
\end{enumerate}
\end{lemma}

\begin{proof}[Proof of Lemma \ref{lem:ell}]
Writing $p_i^{(a)}(\nt)$ to indicate the dependence of the vacancy numbers
on the configuration, one obtains from the algorithm for $\delta$ and
\eqref{eq:vac} that for $1\le a\le n-3$
\begin{equation}\label{eq:vac change}
\begin{split}
p_i^{(a)}(\ntt) = p_i^{(a)}(\nt)
  &-\chi(\ell^{(a-1)}\le i < \ell^{(a)})
   +\chi( \ell^{(a)}\le i < \ell^{(a+1)})\\
  &+\chi(\lb^{(a)}\le i < \lb^{(a-1)})
   -\chi(\lb^{(a+1)} \le i < \lb^{(a)})\\
p_i^{(n-2)}(\ntt) = p_i^{(n-2)}(\nt)
  &-\chi(\ell^{(n-3)} \le i < \ell^{(n-2)})
   +\chi(\ell^{(n-2)}\le i < \min(\ell^{(n-1)},\ell^{(n)}))\\
  &+\chi(\lb^{(n-2)}\le i < \lb^{(n-3)})
   -\chi(\max(\ell^{(n-1)},\ell^{(n)})\le i < \lb^{(n-2)})\\
p_i^{(n-1)}(\ntt) = p_i^{(n-1)}(\nt)
  &-\chi(\ell^{(n-2)}\le i < \ell^{(n-1)})
   +\chi(\ell^{(n-1)}\le i < \lb^{(n-2)})\\
p_i^{(n)}(\ntt) = p_i^{(n)}(\nt)
  &-\chi(\ell^{(n-2)}\le i< \ell^{(n)})
   +\chi(\ell^{(n)}\le i< \lb^{(n-2)}).
\end{split}
\end{equation}

Since $(\nt,\Jt) \in \tj(\RC(\la, B))$, there is a singular string of length
one in the first $k-1$ rigged partitions, so that $\ell^{(a)}=1$ for 
$0\le a \le k-1$ (recall that $\ell^{(0)}=1$ by definition). Similarly, 
$\lt^{(a)}=1$ for $0\le a \le k-2$. Hence $\lt^{(a)}=\ell^{(a+1)}=1$ for 
$0\le a\le k-2$ so that (1) holds in this case.

We proceed by induction on $a$. Fix $k-1\le a\le n-3$.
By induction hypothesis $\lt^{(a-1)}\ge \ell^{(a)} $. By definition
$\lt^{(a)}\ge \lt^{(a-1)}$, so that $\lt^{(a)}\ge \ell^{(a)}$. 
By \eqref{eq:vac change}, $(\ntt,\Jtt)^{(a)}$ has no singular string of length 
$i$ in the range $\ell^{(a)}\le i < \ell^{(a+1)}$ which implies that 
$\lt^{(a)}\ge \ell^{(a+1)}$ for $1\le a \le n-3$. Hence (1) holds.

Next consider $a=n-2$. By construction, $\lt^{(n-2)}\ge \lt^{(n-3)}$ 
and by (1), $\lt^{(n-3)} \ge \ell^{(n-2)}$. Hence $\lt^{(n-2)} \ge \ell^{(n-2)}$.
By \eqref{eq:vac change}, $(\ntt,\Jtt)^{(n-2)}$ has no singular string of length
$i$ in the range $\ell^{(n-2)} \le i < \min(\ell^{(n-1)},\ell^{(n)})$. 
Hence $\lt^{(n-2)} \ge \min(\ell^{(n-1)},\ell^{(n)})$ which proves (2).

Consider $(\ntt,\Jtt)^{(n-1)}$ and $(\ntt,\Jtt)^{(n)}$. By construction,  
$\lt^{(n-1)}\ge \lt^{(n-2)}$ and $\lt^{(n)}\ge \lt^{(n-2)}$.
Using (2) we get $\lt^{(n-1)}\ge \min(\ell^{(n-1)},\ell^{(n)})$ and 
$\lt^{(n)}\ge \min(\ell^{(n-1)},\ell^{(n)})$.
If $\min(\ell^{(n-1)},\ell^{(n)})=\ell^{(n-1)}$, then 
$\lt^{(n-1)}\ge \ell^{(n-1)}$. Now looking at the change in vacancy number 
\eqref{eq:vac change} for $p_i^{(n-1)}$, we conclude that 
$\lt^{(n-1)}\ge \lb^{(n-2)}$. If on the other hand $\min(\ell^{(n-1)},\ell^{(n)})
=\ell^{(n)}$, then $\lt^{(n)}\ge \ell^{(n)}$. Now, by the change in vacancy number
\eqref{eq:vac change} for $p_i^{(n)}$, $\lt^{(n)}\ge \lb^{(n-2)}$. This proves (3).

We prove (4) by downward induction on $a$. Recall that 
$\lbt^{(n-1)}=\max(\lt^{(n-1)},\lt^{(n)})$ by definition. Using (3) we obtain
$\lbt^{(n-1)}\ge \lb^{(n-2)}$. This starts the induction at $a=n-1$. Fix
$2\le a\le n-2$ and assume that $\lbt^{(a+1)}\ge \lb^{(a)}$ by induction hypothesis.
By construction, $\lbt^{(a)}\ge \lbt^{(a+1)}$, so that $\lbt^{(a)}\ge \lb^{(a)}$. 
By the change in vacancy number \eqref{eq:vac change}, $(\ntt,\Jtt)^{(a)}$
has no singular string of length $i$ in the range $\lb^{(a)}\le i < \lb^{(a-1)}$.
This implies that $\lbt^{(a)}\ge \lb^{(a-1)}$ as desired.
\end{proof}

\begin{proof}[Proof of Lemma \ref{lem:s}]
The proof is similar to the proof of Lemma \ref{lem:ell}.
\end{proof}

To show that $\Phi_{\Bt}\circ \tj(\RC(\la,B)) \subset \ts(\Path(\la,B))$, 
let $(\nt,\Jt)\in \tj(\RC(\la,B))$ and $m_k\otimes m_{k-1}m_{k-2}
\cdots m_1 \otimes b'=\Phi_{\Bt}(\nt,\Jt)$. We need to show that either 
$m_k>m_{k-1}$ or $(m_k,m_{k-1})=(n,\nb),(\nb,n)$.

If $m_k=p$ with $p<n$, then $\ell^{(a)}=\infty$ for $a\ge p$ and by 
Lemma \ref{lem:ell} $\tilde{\ell}^{(a)}=\infty$ for $a\ge p-1$. This implies 
that $m_k>m_{k-1}$.

If $m_k=n$, then $\ell^{(n-1)}<\infty$, $\ell^{(n)}=\infty$ and $\lb^{(a)}=\infty$
for $1\le a\le n-2$. By Lemma \ref{lem:ell} this implies $\lt^{(n-1)}=\infty$ and
$\lbt^{(a)}=\infty$ for $1\le a\le n-2$. Hence $m_{k-1}=\overline{n}$ or 
$m_{k-1}<m_k$. The case $m_k=\nb$ is analogous.

If $m_k=\overline{n-1}$, then $\ell^{(n-1)},\ell^{(n)}<\infty$ and $\lb^{(a)}=\infty$
for $1\le a\le n-2$. By Lemma \ref{lem:ell} it follows that either $\lt^{(n-1)}$ or
$\lt^{(n)}$ is $\infty$ so that indeed $m_{k-1}<m_k$.

Finally, if $m_k=\overline{p}$ with $1\le p\le n-2$, then $\lb^{(a)}=\infty$ for 
$1\le a<p$. By Lemma \ref{lem:ell} it follows that $\lbt^{(p)}=\infty$ so that again
$m_k>m_{k-1}$.

Next it needs to be shown that $\Phi_{\Bt}:\tj(\RC(\la,B)) \to \ts(\Path(\la,B))$
is surjective. Let $B=\Bh^{k,1}\otimes B'$ for $2\le k\le n-1$.
Let $b=m_k\otimes m_{k-1} m_{k-2}\cdots m_1\otimes b' \in \ts(\Path(\la,B))$ 
and $(\nt,\Jt)=\Phi^{-1}_{\Bt}(b)$. We need to show that
$(\nt,\Jt)\in \tj(\RC(\la,B))$. In other words, it needs to be shown that
$(\nt,\Jt)^{(a)}$ has a singular string of length one for $1\le a\le k-1$.
By induction on $k$ it follows that $\st^{(a)}=0$ for $0\le a\le k-2$ (remember
that $s^{(0)}=0$ by definition). Hence by Lemma \ref{lem:s}, $s^{(a)}=0$ for 
$0\le a \le k-1$ which implies that there are singular strings of lengths one in 
$(\nt,\Jt)^{(a)}$ for $1\le a\le k-1$.

This completes the proof of case (2).

\subsection*{Proof of \eqref{pt:spinor}}
To prove the existence and uniqueness of a bijection satisfying (3) we
need to show that (a) $\ds$ is well-defined and (b) 
$\Phi_{B'}$ restricts to a bijection $\ds(\RC(\la,B))\to -_s(\Path(\la,B))$.

Let $(\nt,\Jt)\in\RC(\la,B)$ and $(\nt_0,\Jt_0)=\emb_{\RC}(\nt,\Jt)$.
Inductively define $(\nt_j,\Jt_j)=\delta\circ\tj(\nt_{j-1},\Jt_{j-1})$ for
$1\le j\le n$ and let $\ell_j^{(a)}$, $\lb_j^{(a)}$ be the lengths of the
singular strings selected by $\delta$ in $\tj(\nt_j,\Jt_j)$ for $0\le j<n$.
For (a) it needs to be shown that $(\nt_n,\Jt_n)$ is in the image of
$\emb_{\RC}$. This follows from the following lemma:
\begin{lemma}\label{lem:even}
For all $0\le j<n$ we have:
\begin{enumerate}
\item For $1\le a\le n-2$, $\ell_j^{(a)}$ is odd, and 
\begin{equation*}
\ell_j^{(a)} = \begin{cases}
1&\text{if $a<n-j$}\\
\lb_{j-(n-a)}^{(a)}-1 & \text{if $a\ge n-j$.}
\end{cases}
\end{equation*}
In addition for $i=\ell_j^{(a)}-1$, $p_i^{(a)}(\nt_{j+1})$ is even.
\item For $k=n$ and $j$ even, or $k=n-1$ and $j$ odd we have
\begin{equation*}
\begin{split}
&\ell_j^{(n-1)}=\ell_{j-1}^{(n-1)}-1 \quad \text{and $\ell_j^{(n)}$ is even}\\
&\min(\ell_j^{(n-1)},\ell_j^{(n)})=\ell_j^{(n-1)}\\
&\text{$p_i^{(n-1)}(\nt_{j+1})$ is even for $i=\ell_j^{(n-1)}-1$.}
\end{split}
\end{equation*}
Similarly, for $k=n$ and $j$ odd, or $k=n-1$ and $j$ even we have
\begin{equation*}
\begin{split}
&\ell_j^{(n)}=\ell_{j-1}^{(n)}-1 \quad \text{and $\ell_j^{(n-1)}$ is even}\\
&\min(\ell_j^{(n-1)},\ell_j^{(n)})=\ell_j^{(n)}\\
&\text{$p_i^{(n)}(\nt_{j+1})$ is even for $i=\ell_j^{(n)}-1$.}
\end{split}
\end{equation*}
\item For $1\le a\le n-2$, $\lb_j^{(a)}$ is even.
\end{enumerate}
\end{lemma}
Lemma \ref{lem:even} implies that all partitions in $(\nt_n,\Jt_n)$ have only 
parts of even lengths and even riggings so that $(\nt_n,\Jt_n)$ is indeed in the
image of $\emb_{\RC}$. Hence $\ds$ is well-defined.
\begin{proof}[Proof of Lemma \ref{lem:even}]
For concreteness we assume throughout that $k=n$. The case $k=n-1$ is entirely analogous.
For $j=0$, by definition $\tj(\nt_0,\Jt_0)$
has singular strings of lengths 1 in the first $n-1$ rigged partitions. Hence
$\ell_0^{(a)}=1$ for $1\le a\le n-1$. This proves (1). 
Furthermore $\tj(\nt_0,\Jt_0)^{(n)}$ contains only even strings so that 
$\ell_0^{(n)}>\ell_0^{(n-1)}=1$ is even which proves (2). All strings
of length greater than one in $(\nt_0,\Jt_0)^{(a)}$ for $1\le a\le n-2$ are
even so that $\lb_0^{(a)}$ are even which verifies (3).

Suppose the lemma is true for $0\le j'<j$. We will prove it for $j$. Set $b=n-j$.
Clearly $\ell_j^{(a)}=1$ for $1\le a<b$ by the definition of $\tj$. We proceed by
induction on $b\le a<n-2$.

Assume that (1) holds for $p<a$. In particular $\ell_j^{(p)}=\lb_{j-(n-p)}^{(p)}-1$ for 
$b\le p <a$. We want to show that $\ell_j^{(a)}=\lb_{j-(n-a)}^{(a)}-1$. For $a=b$, we have
$\lb_{j-(n-a)}^{(a)}-1\ge 1=\ell_j^{(a-1)}$. Similarly for $a>b$ note that by 
Lemma~\ref{lem:ell} $\lb_{j-(n-a)}^{(a)}-1\ge \lb_{j-(n-a+1)}^{(a-1)}-1 =\ell_j^{(a-1)}$. 
We will show that for $\ell_j^{(a-1)}\le i <\lb_{j-(n-a)}^{(a)}-1$ there is no singular 
string in $(\nt_j,\Jt_j)^{(a)}$. For this range of $i$ the change in vacancy number is
\begin{equation*}
\begin{split}
p_i^{(a)}(\nt_j)=p_i^{(a)}(\nt_0)&-\sum_{k=0}^{j-1}\chi(\ell_k^{(a-1)}\le i<\ell_k^{(a)})
                                 +\sum_{k=0}^{j-1}\chi(\ell_k^{(a)}\le i<\ell_k^{(a+1)})\\ 
                                 &-\sum_{k=0}^{j-1}\chi(\lb_k^{(a+1)}\le i<\lb_k^{(a)})
                                 +\sum_{k=0}^{j-1}\chi(\lb_k^{(a)}\le i<\lb_k^{(a-1)})\\
 =p_i^{(a)}(\nt_0)&+\chi(\ell_{j-1}^{(a)}\le i<\ell_{j-1}^{(a+1)})
                  -\chi(\lb_{a-b}^{(a+1)}\le i<\lb_{a-b}^{(a)}).
\end{split}
\end{equation*}
By induction, $\ell_{j-1}^{(a+1)}= \lb_{a-b}^{(a+1)}-1$, so that the vacancy numbers in
this range change by $\pm 1$ except for $i=\ell_{j-1}^{(a+1)}$.
By induction $\ell_{j-1}^{(a+1)}$ is odd and $\nu_j^{(a)}$ has no odd parts of length
$\lb_{a-b-1}^{(a)}\le i<\lb_{a-b}^{(a)}-1$. In addition the riggings in this range are
even, so that there are no singular strings for $\ell_j^{(a-1)}\le i <\lb_{a-b}^{(a)}-1$
as claimed. Hence  $\ell_j^{(a)}=\lb_{a-b}^{(a)}-1$.

To see that $p_i^{(a)}(\nt_{j+1})$ is even for $i=\ell_j^{(a)}-1$ note that
$\lb_{a-b+m}^{(a-1)}\ge \lb_{a-b}^{(a)}$ and $\lb_{a-b+m+1}^{(a+1)}\ge \lb_{a-b}^{(a)}$
for $m\ge 0$. Furthermore by induction $\ell_{j-m}^{(a-1)}=\lb_{a-b-m-1}^{(a-1)}-1$ for 
$m\ge 0$ and $\ell_{j-m}^{(a+1)}=\lb_{a-b-m+1}^{(a+1)}-1$ for $m\ge 1$, so that all 
parts of length less than $\ell_j^{(a)}$ in $\nu_{j+1}^{(a-1)}$ and $\nu_{j+1}^{(a+1)}$ 
are even. Hence $p_i^{(a)}(\nt_{j+1})$ is even.

For $a=n-2$, we will also show that there is no singular string in 
$(\nt_j,\Jt_j)^{(n-2)}$ of length $\ell_j^{(n-3)}\le i<\lb_{j-2}^{(n-2)}-1$. Similarly 
to before we find that in this range of $i$ the change in vacancy number is
\begin{multline*}
p_i^{(n-2)}(\nt_j)=p_i^{(n-2)}(\nt_0)\\
    +\chi(\ell_{j-1}^{(n-2)}\le i<\min(\ell_{j-1}^{(n-1)},\ell_{j-1}^{(n)}))
    -\chi(\max(\ell_{j-2}^{(n-1)},\ell_{j-2}^{(n)})\le i<\lb_{j-2}^{(n-2)}).
\end{multline*}
If $j$ is even, by induction $\min(\ell_{j-1}^{(n-1)},\ell_{j-1}^{(n)})
=\ell_{j-1}^{(n)}$, $\max(\ell_{j-2}^{(n-1)},\ell_{j-2}^{(n)})=\ell_{j-2}^{(n)}$
and $\ell_{j-1}^{(n)}=\ell_{j-2}^{(n)}-1$ with $\ell_{j-2}^{(n)}$ is even, so that
the vacancy number changes by $\pm 1$ except for $i=\ell_{j-2}^{(n)}-1$. However,
there is no odd string in $\nu_j^{(n-2)}$ for $\lb_{j-3}^{(n-3)}\le i<
\lb_{j-2}^{(n-2)}-1$ and all riggings in this range are even. Hence
$\ell_j^{(n-2)}=\lb_{j-2}^{(n-2)}-1$ as desired. The case $j$ odd is analogous.

To verify that $p_i^{(n-2)}(\nt_{j+1})$ is even for $i=\ell_j^{(n-2)}-1$
let us again assume that $j$ is even. Note that $\ell_{j-1}^{(n-1)}\ge \lb_{j-2}^{(n-2)}$ 
by Lemma \ref{lem:ell} and by induction $\ell_{j-k+1}^{(n)}=\ell_{j-k}^{(n)}-1$ for
$k>0$ even and $\ell_{j-k+1}^{(n-1)}=\ell_{j-k}^{(n-1)}-1$ for $k>1$ odd. Hence
$\nu_{j+1}^{(n-1)}$ and $\nu_{j+1}^{(n)}$ contain only even strings of length
less than $\ell_j^{(n-2)}$. Similarly, $\nu_{j+1}^{(n-3)}$ contains only even
strings of length less than $\ell_j^{(n-2)}$ which proves the assertion. The case
$j$ odd is again analogous. This concludes the proof of (1).

For the proof of (2) assume that $j$ is even. To prove 
$\ell_j^{(n-1)}=\ell_{j-1}^{(n-1)}-1$ we show that there are no singular strings of 
length $\ell_j^{(n-2)}\le i<\ell_{j-1}^{(n-1)}-1$ in $(\nt_j,\Jt_j)^{(n-1)}$. By 
induction $\min(\ell_{j-2}^{(n-1)},\ell_{j-2}^{(n)})=\ell_{j-2}^{(n-1)}$, so that by 
Lemma \ref{lem:ell} $\ell_{j-1}^{(n-1)}\ge \lb_{j-2}^{(n-2)}=\ell_j^{(n-2)}+1$. 
Furthermore, $\ell_j^{(n-2)}\ge \ell_{j-1}^{(n-2)}$. Hence for
$\ell_j^{(n-2)}\le i<\ell_{j-1}^{(n-1)}-1$ the vacancy number changes by 
\begin{equation*}
\begin{split}
p_i^{(n-1)}(\nt_j)=& p_i^{(n-1)}(\nt_0)
                   -\sum_{k=0}^{j-1}\chi(\ell_k^{(n-2)}\le i<\ell_k^{(n-1)})
                   +\sum_{k=0}^{j-1}\chi(\ell_k^{(n-1)}\le i<\lb_k^{(n-2)})\\
                  =& p_i^{(n-1)}(\nt_0)-\chi(\ell_{j-1}^{(n-2)}\le i<\ell_{j-1}^{(n-1)}).
\end{split}
\end{equation*}
Since $p_i^{(n-1)}(\nt_0)$ is even and all riggings in this range are even, this implies
that there are no singular strings. Hence $\ell_j^{(n-1)}=\ell_{j-1}^{(n-1)}-1$
as claimed.

By induction $\min(\ell_{j-1}^{(n-1)},\ell_{j-1}^{(n)})=\ell_{j-1}^{(n)}$, so that
by Lemma \ref{lem:ell} $\ell_j^{(n)}\ge \lb_{j-1}^{(n-2)}$. On the other hand
$\ell_j^{(n-1)}=\ell_{j-1}^{(n-1)}-1<\lb_{j-1}^{(n-2)}$ which proves that
$\min(\ell_j^{(n-1)},\ell_j^{(n)})=\ell_j^{(n-1)}$. 
Note that $\ell_{j'}^{(n)}\le \lb_{j-1}^{(n-2)}$ for $j'<j$ so that there are only
even strings of length greater or equal to $\lb_{j-1}^{(n-2)}$ in $\nu_j^{(n)}$.
Hence $\ell_j^{(n)}$ must be even.
All strings of length less than $\ell_j^{(n-1)}$ in $\nu_{j+1}^{(n-2)}$ are even
by induction. Hence $p_i^{(n-1)}(\nt_{j+1})$ is even for $i=\ell_j^{(n-1)}-1$.

The case $j$ is odd is analogous. This concludes the proof of (2).

For (3) note that $\lb_j^{(a)}\ge \lb_{j-1}^{(a-1)}$ by Lemma \ref{lem:ell}.
There are only even strings of length greater or equal to $\lb_{j-1}^{(a-1)}$
in $\nu_j^{(a)}$ which proves that $\lb_j^{(a)}$ is even.
\end{proof}

To prove (b) we need to show that the crystal element defined by $\ds$ is indeed
an element of $B^{n,1}$ (resp. $B^{n-1,1}$) and that $\ds^{-1}$ exists and is
well-defined. The proof of the latter is similar to the proof of Lemma \ref{lem:even}
and is omitted here. To show that the step defined by $\ds$ is in $B^{n,1}$
(resp. $B^{n-1,1}$) we need to show that \eqref{eq:spinor} holds. 
For concreteness let us assume that the crystal is $B^{n,1}$.
The condition $m_h<m_{h+1}$ follows from (2) of Theorem \ref{thm:bij} by the
definition of $\ds$.

Let $j$ be maximal such that both $\ell_{j-1}^{(n-1)},\ell_{j-1}^{(n)}<\infty$.
If $j$ is even, $\ell_j^{(n-1)}=\ell_{j-1}^{(n-1)}-1<\infty$ by Lemma \ref{lem:even},
so that $\ell_j^{(n)}=\infty$. This implies that $m_{n-j}=n$. Similarly, if
$j$ is odd, $\ell_j^{(n)}=\ell_{j-1}^{(n)}-1<\infty$ and
$\ell_j^{(n-1)}=\infty$, so that $m_{n-j}=\overline{n}$. This proves
point 3 of \eqref{eq:spinor}. It also shows that either $n$ or $\overline{n}$ occurs
in $b=m_n\ldots m_1$. To see that not both $n$ and $\overline{n}$ can occur, we note 
that if $j$ is even then $\ell_{j+1}^{(n)}=\ell_j^{(n)}-1=\infty$ and
$\min(\ell_{j+1}^{(n-1)},\ell_{j+1}^{(n)})=\ell_{j+1}^{(n)}$ by Lemma \ref{lem:even},
so that $\ell_{j+1}^{(n-1)}=\infty$ as well. Hence $m_{n-j-1}\le n-1$. The case
$j$ odd is analogous. 

It remains to show that precisely one of $p$ and $\overline{p}$ occurs in $b$ for 
$1\le p\le n-1$. Suppose $m_{n-j}=\overline{p}$. Then $\lb_j^{(p-1)}=\infty$ and
$\lb_j^{(p)}<\infty$. By Lemma \ref{lem:ell} $\lb_{j-1}^{(p-1)}\le \lb_j^{(p)}<\infty$
and $\lb_{j+1}^{(p)}\ge \lb_j^{(p-1)}=\infty$. Setting $a=n-p+j$ it follows from 
Lemma \ref{lem:even} that
\begin{align*}
\ell_a^{(p-1)}&=\lb_{j-1}^{(p-1)}-1<\infty & \ell_a^{(p)}&=\lb_j^{(p)}-1<\infty\\
\ell_{a+1}^{(p-1)}&=\lb_j^{(p-1)}-1=\infty & \ell_{a+1}^{(p)}&=\lb_{j+1}^{(p)}-1=\infty.
\end{align*} 
Hence $p$ does not occur in $b$.

Now suppose $\overline{p}\not\in b$. We show that then $p\in b$. 
If $\overline{p}\not\in b$ there must be some $j$ such that 
$\lb_{j-1}^{(p)},\lb_{j-1}^{(p-1)}<\infty$ and $\lb_j^{(p)},\lb_j^{(p-1)}=\infty$.
Setting as before $a=n-p+j$ this implies by Lemma \ref{lem:even} that
$\ell_a^{(p-1)}=\lb_{j-1}^{(p-1)}-1<\infty$ and $\ell_a^{(p)}=\lb_j^{(p)}-1=\infty$.
Hence $m_{n-a}=p$. This proves point 2 of \eqref{eq:spinor}.

This concludes the proof of Theorem \ref{thm:bij}.

\section{Proof of Theorem \ref{thm:emb}}
\label{app:emb proof}

In the following let us suppress the weight $\la$ and write $\RC(B)$ for
$\RC(\la,B)$. It is straightforward to check that the diagram
\begin{equation}\label{eq:commute RC}
\xymatrix{
  {\RC(B^{k,1}\otimes B)} \ar[rr]^{\emb_{RC}} \ar[ddd]_{\delta\circ\tj}
 &&{\RC(B^{k,1}\otimes B^{k,1}\otimes \Bh)} \ar[d]^{\delta\circ\tj}\\
 &&{\RC(B^{k-1,1}\otimes B^{k,1}\otimes \Bh)}\ar[d]^{\id}\\
 &&{\RC(B^{k,1}\otimes B^{k-1,1}\otimes \Bh)}\ar[d]^{\delta\circ\tj}\\
  {\RC(B^{k-1,1}\otimes B)} \ar[rr]^{\emb_{RC}}
 &&{\RC(B^{k-1,1}\otimes B^{k-1,1}\otimes \Bh)}
}
\end{equation}
commutes.

We will show here that the crystal analogon to this diagram commutes:
\begin{equation}\label{eq:commute B}
\xymatrix{
  {B^{k,1}} \ar[rr]^{\emb_B} \ar[ddd]_{-\circ\ts}
 &&{B^{k,1}\otimes B^{k,1}} \ar[d]^{-\circ\ts}\\
 &&{B^{k-1,1}\otimes B^{k,1}}\ar[d]^{R}\\
 &&{B^{k,1}\otimes B^{k-1,1}}\ar[d]^{-\circ\ts}\\
  {B^{k-1,1}} \ar[rr]^{\emb_B}
 &&{B^{k-1,1}\otimes B^{k-1,1}.}
}
\end{equation}
By the definition of $\Phi$ and Proposition \ref{prop:corresp}, the commutative 
diagrams \eqref{eq:commute RC} and \eqref{eq:commute B} imply Theorem \ref{thm:emb}.

First we verify \eqref{eq:commute B} for classically highest weight elements.
All classically highest weight elements $b\in B^{k,1}$ are of the form
$b=k\cdots 21$ or $b=\overline{\ell}\cdots \overline{p}p\cdots 21$ for some 
$1\le \ell\le p$ such that $2p-\ell+1=k$. For $b=k\cdots 21$ everything is
determined by weight considerations and it is easy to see that both paths
in the diagram \eqref{eq:commute B} yield $(k-1)\cdots 21\otimes (k-1)\cdots 21$.
For $b=\overline{\ell}\cdots \overline{p}p\cdots 21$ explicit calculations show that
$\emb_B(b)=\overline{\ell}\cdots\overline{k}(\ell-1)\cdots 21 \otimes k \cdots 21$.
Hence 
\begin{equation*}
\begin{split}
&-\circ\ts(b)= \overline{\ell+1}\cdots \overline{p}\;p\cdots 21\\
&\emb_B\circ-\circ\ts(b)=\overline{\ell+1}\cdots \overline{k-1}\;\ell\cdots 21
 \otimes (k-1)\cdots 21\\
&-\circ\ts\circ\emb_B(b)=\overline{\ell+1}\cdots \overline{k}\;(\ell-1)\cdots 21
 \otimes k\cdots 21.
\end{split}
\end{equation*}
{}From its definition, that $R$ commutes with all crystal operators, one may
check that
\begin{equation*}
R(\overline{\ell+1}\cdots \overline{k}\;(\ell-1)\cdots 21 \otimes k\cdots 21)
= \overline{\ell}\cdots \overline{k-1}\; \ell\cdots 21\otimes (k-1)\cdots 21.
\end{equation*}
Hence indeed $\emb_B\circ-\circ\ts(b)=-\circ\ts\circ R\circ -\circ\ts\circ\emb_B(b)$.

Next we show that if the diagram commutes for $b\in B^{k,1}$ then it commutes
for $f_ib$ for $1\le i\le n$ provided that $f_ib\neq 0$. 

Suppose that $b=\beta b'$ where $\beta$ is a letter in the alphabet 
$\A=\{1,2,\ldots,n,\overline{n},\ldots,\overline{2},\overline{1}\}$ and $b'$ is a column 
of height $k-1$. We claim that
\begin{align}
\label{eq:1}
\emb_B(b) &= \beta a'\otimes a\\
\label{eq:2}
R(a'\otimes a) &= \beta \tilde{a} \otimes \tilde{a}'
\end{align}
for some column $a'$ of height $k-1$, $a$ of height $k$, and columns
$\tilde{a}$ and $\tilde{a}'$ of height $k-1$.
Equations \eqref{eq:1} and \eqref{eq:2} are certainly true for classically highest
weight elements $b$ as is clear from the previous explicit calculations.
Now suppose \eqref{eq:1} and \eqref{eq:2} hold for $b$. Then they hold
for $f_ib$ for $1\le i\le n$ if $f_ib\neq 0$ since $\emb_B(B^{k,1})$ is aligned and
$R$ is a crystal isomorphism so that $f_i R=R f_i$.

If $f_i$ does not act on $\beta$ in $b$, then 
$\emb_B\circ-\circ\ts(f_ib)=f_i^2\circ\emb_B\circ-\circ\ts(b)$ and
\begin{equation*}
\begin{split}
 &-\circ\ts\circ R\circ-\circ\ts\circ\emb_B(f_ib)\\
=&-\circ\ts\circ R\circ-\circ\ts\circ f_i^2\circ\emb_B(b)\\
=&-\circ\ts\circ R\circ f_i^2\circ-\circ\ts\circ\emb_B(b)\\
=&-\circ\ts\circ f_i^2 \circ R\circ-\circ\ts\circ\emb_B(b)\\
=&f_i^2\circ-\circ\ts\circ R\circ-\circ\ts\circ\emb_B(b)
\end{split}
\end{equation*}
by \eqref{eq:1} and \eqref{eq:2} so that the diagram still commutes.
Similarly, if $f_i$ acts on $\beta$ in $b$, then 
$\emb_B\circ-\circ\ts(f_ib)=\emb_B\circ-\circ\ts(b)$ and also 
\begin{equation*}
\begin{split}
 &-\circ\ts\circ R\circ-\circ\ts\circ\emb_B(f_ib)\\
=&-\circ\ts\circ R\circ-\circ\ts\circ f_i^2\circ\emb_B(b)\\
=&-\circ\ts\circ R\circ f_i\circ-\circ\ts\circ\emb_B(b)\\
=&-\circ\ts\circ f_i \circ R\circ-\circ\ts\circ\emb_B(b)\\
=&-\circ\ts\circ R\circ-\circ\ts\circ\emb_B(b)
\end{split}
\end{equation*}
by \eqref{eq:1} and \eqref{eq:2} so that the diagram still commutes.
This concludes the proof of Theorem \ref{thm:emb}.

\section{Proof of Lemma \ref{lem:delta-delta}}
\label{app:delta-delta}

One may easily verify that (see also \cite[Eq. (3.10)]{KSS:2002})
\begin{equation}\label{eq:p ineq}
\begin{split}
-p_{i-1}^{(a)}+2p_i^{(a)}-p_{i+1}^{(a)}
=&-\sum_{b\in J} \left( \alpha_a \mid \alpha_b\right) m_i^{(b)}+\delta_{i,1} L^a\\
\ge & -\sum_{b\in J} \left( \alpha_a \mid \alpha_b\right) m_i^{(b)},
\end{split}
\end{equation}

The proof of $[\delta,\dt]=0$ follows from Lemmas \ref{lem:ell s}
and \ref{lem:J} below. The proof of these lemmas is very technical and follows
similar arguments as in~\cite[Appendix A]{KSS:2002}. Details are available upon
request.

Let $(\nt,\Jt)\in \RC(\lambda,B)$ where $B=(B^{1,1})^{\otimes 2} \otimes B'$.
The following notation is used:
\begin{equation*}
\begin{split}
\delta(\nt,\Jt)&=(\nud,\Jd)\\
\dt(\nt,\Jt)&=(\nut,\Jtil)\\
\dt\circ \delta(\nt,\Jt)&=(\nudt,\Jdt)\\
\delta\circ \dt(\nt,\Jt)&=(\nutd,\Jtd).
\end{split}
\end{equation*}

Furthermore, let $\{\ld^{(k)},\ldb^{(k)}\}$, $\{\lt^{(k)},\ltb^{(k)}\}$,
$\{\ldt^{(k)},\ldtb^{(k)}\}$ and $\{\ltd^{(k)},\ltdb^{(k)}\}$
be the lengths of the strings that are shortened in the transformations 
$(\nt,\Jt)\mapsto (\nud,\Jd)$, $(\nt,\Jt)\mapsto (\nut,\Jtil)$,
$(\nud,\Jd)\mapsto (\nudt,\Jdt)$ and $(\nut,\Jtil)\mapsto (\nutd,\Jtd)$,
respectively. (In Section~\ref{sec:delta} the lengths of the $\s$-strings was denoted
by $\lb$).
We call the strings, whose lengths are labeled by an $\ell$,
$\ell$-strings and those labeled by an $\s$, $\s$-strings.

\begin{lemma}\label{lem:ell s}
The following cases occur at $(\nt,\Jt)^{(k)}$:
\begin{itemize}
\item[I.] \textbf{Nontwisted case.} In this case the $\ell$-string selected by
$\delta$ (resp. $\dt$) in $(\nt,\Jt)^{(k)}$ is different from the $\s$-string selected 
by $\dt$ (resp. $\delta$) in $(\nt,\Jt)^{(k)}$.
For the $\ell$-strings one of the following must hold:
\begin{itemize}
\item[($\ell$a)] \textbf{Generic case.} If $\delta$ and $\dt$ do not select the same $\ell$-string,
then $\ltd^{(k)}=\ld^{(k)}$ and $\ldt^{(k)}=\lt^{(k)}$.
\item[($\ell$b)] \textbf{Doubly singular case.} In this case $\delta$ and $\dt$ select the same
$\ell$-string, so that $\ld^{(k)}=\lt^{(k)}=:\ell$. Then
\begin{itemize}
\item[(1)] If $\ldt^{(k)}<\ell$ (or $\ltd^{(k)}<\ell$) then $\ldt^{(k)}=\ltd^{(k)}=\ell-1$
and $m_{\ell-1}^{(k+1)}=0$ for $k<n-2$, $m_{\ell-1}^{(n-1)}=m_{\ell-1}^{(n)}=0$ for $k=n-2$
and $m_{\ell-1}^{(n-2)}=0$ for $k=n-1,n$.
\item[(2)] If $\ldt^{(k)}=\ell$ (or $\ltd^{(k)}=\ell$) then case I.($\ell\s$)(1')
(or I.($\ell\s$)(1)) holds or $\ldt^{(k)}=\ltd^{(k)}=\ell$.
\item[(3)] If $\ldt^{(k)}>\ell$ (or $\ltd^{(k)}>\ell$) then case I.($\ell s$)(1',2)
(or I.($\ell\s$)(1,2)) holds or
$\ldt^{(k)}=\ltd^{(k)}$ and $\ldt^{(k)}\le \lt^{(k+1)}$, $\ltd^{(k)}\le \ld^{(k+1)}$ for $k<n-2$,
$\ldt^{(n-2)}\le \min\{\lt^{(n-1)},\lt^{(n)}\}$, $\ltd^{(n-2)}\le \min\{\ld^{(n-1)},\ld^{(n)}\}$ 
for $k=n-2$, and $\ldt^{(k)}\le \ltb^{(n-2)}$, $\ltd^{(k)}\le \ldb^{(n-2)}$ for $k=n-1,n$.
\end{itemize}
\end{itemize}
For the $\s$-strings, case I.($\ell s$) holds or one the following must hold:
\begin{itemize}
\item[($\s$a)] \textbf{Generic case.} If $\delta$ and $\dt$ do not select the same $\s$-string,
then $\ltdb^{(k)}=\ldb^{(k)}$ and $\ldtb^{(k)}=\ltb^{(k)}$.
\item[($\s$b)] \textbf{Doubly singular case.} In this case $\delta$ and $\dt$ select the same 
$\s$-string, so that $\ldb^{(k)}=\ltb^{(k)}=:\s$. Then
\begin{itemize}
\item[(1)] If $\ldtb^{(k)}<\s$ (or $\ltdb^{(k)}<\s$) then $\ldtb^{(k)}=\ltdb^{(k)}=\s-1$
and $m_{\s-1}^{(k-1)}=0$.
\item[(2)] If $\ldtb^{(k)}=\s$ (or $\ltdb^{(k)}=\s$) then $\ldtb^{(k)}=\ltdb^{(k)}=\s$.
\item[(3)] If $\ldtb^{(k)}>\s$ (or $\ltdb^{(k)}>\s$) then $\ldtb^{(k)}=\ltdb^{(k)}$,
$\ldtb^{(k)}\le \ltb^{(k-1)}$ and $\ltdb^{(k)}\le \ldb^{(k-1)}$.
\end{itemize}
\item[($\ell s$)] \textbf{Mixed case.} One of the following holds:
\begin{itemize}
\item[(1)] $\ld^{(k)}=\lt^{(k)}=:\ell$, $\ltd^{(k)}=\ldb^{(k)}=\ld^{(k+1)}=:\ell'$,
$\ldt^{(k)}=\ltdb^{(k)}=:\ell''$, $\ltb^{(k)}=\ldtb^{(k)}$ or possibly the same conditions for
$\ell$ and $\ell'$, $\ldt^{(k)}=\ltb^{(k)}=\ltb^{(k+1)}=\ell''$, 
$\ltdb^{(k)}=\ldtb^{(k)}=\ell'''$, $m_{\ell''}^{(k-1)}=0$, $m_{\ell''}^{(k)}=1$, 
$m_{\ell''}^{(k+1)}=2$ if case I.($\ell\s$)(1) does not hold at $k-1$. 
Furthermore, either $\ldt^{(k)}\le \lt^{(k+1)}$ or case I.($\ell s$)(1) holds at $k+1$ 
with the same values of $\ell'$ and $\ell''$.
Similarly, either $\ltdb^{(k)}\le \ldb^{(k-1)}$ or case I.($\ell s$)(1) holds at $k-1$
with the same values of $\ell'$ and $\ell''$.
\item[(1')] $\ld^{(k)}=\lt^{(k)}=:\ell$, $\ldt^{(k)}=\ltb^{(k)}=\lt^{(k+1)}=:\ell'$,
$\ltd^{(k)}=\ldtb^{(k)}=:\ell''$, $\ldb^{(k)}=\ltdb^{(k)}$ or possibly the same conditions for
$\ell$ and $\ell'$, $\ltd^{(k)}=\ldb^{(k)}=\ldb^{(k+1)}=\ell''$, 
$\ldtb^{(k)}=\ltdb^{(k)}=\ell'''$, $m_{\ell''}^{(k-1)}=0$, $m_{\ell''}^{(k)}=1$, 
$m_{\ell''}^{(k+1)}=2$ if case I.($\ell\s$)(1') does not hold at $k-1$. 
Furthermore, either $\ltd^{(k)}\le \ld^{(k+1)}$ or case I.($\ell s$)(1') holds at $k+1$ 
with the same values of $\ell'$ and $\ell''$.
Similarly, either $\ldtb^{(k)}\le \ltb^{(k-1)}$ or case I.($\ell s$)(1') holds at $k-1$
with the same values of $\ell'$ and $\ell''$.
\item[(2)] For $k<n-2$ (resp. $k=n-2$) $\ld^{(k)}=\lt^{(k)}=:\ell$, 
$\ldb^{(k)}=\ltb^{(k)}=\ldb^{(k+1)}=\ltb^{(k+1)}=:\ell'$
(resp. $\ldb^{(k)}=\ltb^{(k)}=\ld^{(n-1)}=\ld^{(n)}=\lt^{(n-1)}=\lt^{(n)}=\ell'$),
$\ltd^{(k)}=\ldt^{(k)}=\ell''$, $\ltdb^{(k)}=\ldtb^{(k)}:=\ell'''$ and case I.($\ell s$)(2) 
holds at $k+1$ (resp. $n-1$ and $n$) with the same values of $\ell'$ and $\ell''$ and 
$\ell=\ell'$, $\ell'''=\ell''$. 
Also, either $\ldtb^{(k)}\le \ltb^{(k-1)}$ and $\ltdb^{(k)}\le \ldb^{(k-1)}$ or 
case I.($\ell s$)(2) holds at $k-1$ with the same values of $\ell'$ and $\ell''$.
For $k=n-1,n$, $\ld^{(k)}=\lt^{(k)}=\ldb^{(n-2)}=\ltb^{(n-2)}=\ell'$ and 
$\ldt^{(k)}=\ltd^{(k)}=\ell''$. In addition case I.($\ell\s$)(2) holds at $n-2$
with the same values of $\ell'$ and $\ell''$.
\end{itemize}
\end{itemize}
\item[II.] \textbf{Twisted case.} In this case the $\ell$-string in $(\nt,\Jt)^{(k)}$ selected 
by $\delta$ is the same as the $\s$-string selected by $\dt$ or vice versa.
In the first case $\ld^{(k)}=\ltb^{(k)}=:\ell$. Then $\lt^{(k)}=\ldt^{(k)}$ and one of the 
following holds:
\begin{itemize}
\item[(1)] If $\ltd^{(k)}<\ell$, then $\ltd^{(k)}=\ldtb^{(k)}=\ell-1$,
$m_{\ell-1}^{(k+1)}=0$ or $m_{\ell-1}^{(k+1)}(\nut)=0$, and $m_{\ell-1}^{(k-1)}=0$ or
$m_{\ell-1}^{(k-1)}(\nud)=0$. Furthermore $\ldb^{(k)}=\ltdb^{(k)}$.
\item[(2)] If $\ltd^{(k)}=\ell$, then $\ltd^{(k)}=\ldtb^{(k)}=\ell$ and 
$\ldb^{(k)}=\ltdb^{(k)}$.
\item[(3)] If $\ltd^{(k)}>\ell$, then
\begin{itemize}
\item[(i)] $\ltd^{(k)}=\ldtb^{(k)}$ and $\ldb^{(k)}=\ltdb^{(k)}$, or
\item[(ii)] $\ltd^{(k)}=\ldb^{(k)}$ and $\ldtb^{(k)}=\ltdb^{(k)}\le \ldb^{(k-1)}$.
\end{itemize}
Furthermore, either $\ltd^{(k)}\le \ld^{(k+1)}$ or $\ld^{(k)}=\ld^{(k+1)}$, 
$\ltd^{(k)}=\ltd^{(k+1)}$, $m_\ell^{(k+1)}=1$ and Case II.(3)(i) holds at $k+1$.
Similarly, either $\ldtb^{(k)}\le \ltb^{(k-1)}$ or $\ld^{(k)}=\ld^{(k-1)}$, 
$\ltd^{(k)}=\ltd^{(k-1)}$, $m_\ell^{(k-1)}=1$ and Case II.(3) holds at $k-1$.
\end{itemize}
If the $\ell$-string in $(\nt,\Jt)^{(k)}$ selected by $\dt$ is the same as the 
$\s$-string selected by $\delta$, then $\lt^{(k)}=\ldb^{(k)}=:\ell$. In this 
case $\ld^{(k)}=\ltd^{(k)}$ and one of the following holds:
\begin{itemize}
\item[(1')] If $\ldt^{(k)}<\ell$, then $\ldt^{(k)}=\ltdb^{(k)}=\ell-1$,
$m_{\ell-1}^{(k+1)}=0$ or $m_{\ell-1}^{(k+1)}(\nud)=0$, and $m_{\ell-1}^{(k-1)}=0$ or
$m_{\ell-1}^{(k-1)}(\nut)=0$. Furthermore $\ltb^{(k)}=\ldtb^{(k)}$.
\item[(2')] If $\ldt^{(k)}=\ell$, then $\ldt^{(k)}=\ltdb^{(k)}=\ell$ and $\ltb^{(k)}=\ldtb^{(k)}$.
\item[(3')] If $\ldt^{(k)}>\ell$, then
\begin{itemize}
\item[(i)] $\ldt^{(k)}=\ltdb^{(k)}$ and $\ltb^{(k)}=\ldtb^{(k)}$, or
\item[(ii)] $\ldt^{(k)}=\ltb^{(k)}$ and $\ltdb^{(k)}=\ldtb^{(k)}\le \lt^{(k-1)}$.
\end{itemize}
Furthermore, either $\ldt^{(k)}\le \lt^{(k+1)}$ or $\lt^{(k)}=\lt^{(k+1)}$, 
$\ldt^{(k)}=\ldt^{(k+1)}$, $m_\ell^{(k+1)}=1$ and Case II.(3')(i) holds at $k+1$.
Similarly, either $\ltdb^{(k)}\le \ldb^{(k-1)}$ or $\lt^{(k)}=\lt^{(k-1)}$, 
$\ldt^{(k)}=\ldt^{(k-1)}$, $m_\ell^{(k-1)}=1$ and Case II.(3') holds at $k-1$.
\end{itemize}
\end{itemize}
\end{lemma}

\begin{lemma}\label{lem:J}
$\Jdt=\Jtd$.
\end{lemma}

\section{Combinatorial $R$-matrix on highest weight elements}
\label{app:R matrix}

Here we give the combinatorial $R$-matrix 
$R:B^{k',1}\otimes B^{k,1} \to B^{k,1}\otimes B^{k',1}$
explicitly for highest weight elements. 
The following notation is employed. A column with entries
$\cdots a_4\cdots a_3 a_2 \cdots a_1$, where $a_2\cdots a_1$ denotes a
sequence of consecutive integers $a_2\, a_2-1\, a_2-2\,\ldots a_1+1\, a_1$, is compactly 
written as $\left[a_1..a_2 \mid a_3..a_4 \mid \cdots \right]$.
We use the convention that a segment $[..\mid a..b \mid ..]$ is empty if 
$b<a$.

Since $R$ is an involution, it suffices to state $R$ for the case $k'\le k$.
For $k'=k$ we have $R=\mathrm{id}$.
Hence for $R:B^{n,1}\otimes B^{n,1} \to B^{n,1}\otimes B^{n,1}$ the $R$-matrix
is the identity. The highest weight elements are given by 
$[1..n-2\ell \mid \ol{n}..\ol{n-2\ell+1}]\otimes [1..n]$ for $0\le \ell\le n/2$
with local energy $H=\ell$. The same results hold for 
$R:B^{n-1,1}\otimes B^{n-1,1}\to B^{n-1,1} \otimes B^{n-1,1}$ with $n$ and 
$\ol{n}$ interchanged.

For $R:B^{n-1,1}\otimes B^{n,1} \to B^{n,1} \otimes B^{n-1,1}$
the highest weight elements get mapped as follows:
\begin{equation*}
\begin{split}
&[1..n-2\ell-1 \mid \ol{n}..\ol{n-2\ell}] \otimes [1..n]\\
\overset{R}{\longleftrightarrow}&
[1..n-2\ell-1 \mid n \mid \ol{n-1}..\ol{n-2\ell}] \otimes [1..n-1 \mid \ol{n}]
\end{split}
\end{equation*}
where $0\le \ell<n/2$. The local energy is $H=\ell$.

Next consider the case $R:B^{k',1} \otimes B^{n,1} \to B^{n,1} \otimes B^{k',1}$
where $k'\le n-2$. The highest weight elements are 
$[1..\ell \mid \ol{n}n...\ol{n}n \mid \ol{\ell}..\ol{r}] \otimes [1..n]$
and $[1..\ell \mid \ol{n}n...\ol{n} \mid \ol{\ell}..\ol{r}] \otimes [1..n]$
where $r\le \ell+1$ and the pieces $\ol{n}n..\ol{n}n$ and $\ol{n}n..\ol{n}$ could 
have length zero. Then
\begin{equation*}
\begin{split}
[1..\ell \mid \ol{n}n...\ol{n}n \mid \ol{\ell}..\ol{r}] \otimes [1..n]
\overset{R}{\longleftrightarrow}&
[1..r-1 \mid \ell+d+1..n \mid \ol{\ell+d}..\ol{r}]\\
\otimes &\left[1..\frac{k'+\ell+d}{2} \mid \ol{\frac{k'+\ell+d}{2}}..\ol{\ell+d+1}\right]
\end{split}
\end{equation*}
with $H=\frac{k'-\ell+d}{2}$ and 
\begin{equation*}
\begin{split}
[1..\ell \mid \ol{n}n...\ol{n} \mid \ol{\ell}..\ol{r}] \otimes [1..n]
\overset{R}{\longleftrightarrow}&
[1..r-1 \mid \ell+d+2..n-1 \mid \ol{n} \mid \ol{\ell+d+1}..\ol{r}]\\
\otimes 
&\left[1..\frac{k'+\ell+d+1}{2} \mid \ol{\frac{k'+\ell+d+1}{2}}..\ol{\ell+d+2}\right]
\end{split}
\end{equation*}
with $H=\frac{k'-\ell+d+1}{2}$.

Similarly, the case $R:B^{k',1} \otimes B^{n-1,1} \to B^{n-1,1} \otimes B^{k',1}$
with $k'\le n-2$ is obtained from the above by interchanging $n$ and $\ol{n}$ everywhere.

The remaining cases are all of the form $R:B^{k',1}\otimes B^{k,1} \to
B^{k,1} \otimes B^{k',1}$ with $1\le k'\le k\le n-2$. The highest weight elements
$u'\otimes u\in B^{k',1}\otimes B^{k,1}$ are summarized in the following table:
\newline
\begin{tabular}{|c|l|l|}
\hline
Case & $u'\otimes u$ & parameter ranges\\
\hline &&\\[-2mm]
1 & $\left[1..\ell \mid k+1..p \mid \ol{p}..\ol{q} \mid \ol{\ell}..\ol{r}\right]
\otimes \left[1..k\right]$
& $\ell\ge 0$; $p\ge k$; $q\le p+1$; $r\le \ell+1$\\[2mm]
2a & $[1..\ell \mid \ol{\ell}..\ol{r}] \otimes [1..a \mid \ol{a}..\ol{b}]$
& $a\ge b$; $\ell\ge 1$; $r\le \ell+1$\\[2mm]
2b & $[1..\ell \mid \ol{\ell}..\ol{r} \mid \ol{b-1}..\ol{q}] 
 \otimes [1..a \mid \ol{a}..\ol{b}]$
& $a\ge b$; $\ell\ge r>b>q\ge 1$\\[2mm]
3 & $[1..\ell \mid b..p \mid \ol{\ell}..\ol{r}] \otimes [1..a\mid \ol{a}..\ol{b}]$
& $a\ge b$; $\ell\ge 0$; $p\ge b>\ell+1\ge r$\\[2mm]
4 & $[1..\ell \mid b..p \mid \ol{p}..\ol{q} \mid \ol{\ell}..\ol{r}] \otimes 
 [1..a \mid \ol{a}..\ol{b}]$
& $a\ge b$; $\ell\ge 0$; $p\ge b,q$; $r\le \ell+1$\\[2mm]
5a & $[1..\ell \mid \ol{b-1}..\ol{q}] \otimes [1..a \mid \ol{a}..\ol{b}]$
& $a\ge b>q$; $\ell\ge 0$\\[2mm]
5b & $[1..\ell \mid \ol{b-1}..\ol{q} \mid \ol{\ell}..\ol{r}] 
 \otimes [1..a \mid \ol{a}..\ol{b}]$
& $a\ge b>q>\ell+1$; $\ell \ge r\ge 1$\\[2mm]
\hline
\end{tabular}
\newline
In cases 1 and 4 $..p\mid \ol{p}..$ may also be replaced by
$..n\ol{n}\cdots n\ol{n}..$ or $..\ol{n}n\cdots \ol{n}n..$.

The answers for $R(u'\otimes u)$ in the various cases are given in the tables below.
For brevity, we use the definition $d=\ell-r+1$. By weight considerations 
$R([1..k']\otimes [1..k])=[1..k]\otimes [1..k']$.
All answers are proven by induction.
As an example, we demonstrate how to derive Case 2a with $r\le \ell$ and $k'<b\le a-d$ 
by induction on $d$. We assume that Case 2a holds for $d=0$ (which can also be proven by
induction).

First consider $d=1$ so that $u=[1..\ell \mid \ol{\ell}]\otimes [1..a\mid \ol{a}..\ol{b}]$.
Then $k'<b$ implies that $\ell<b-1$ so that
\begin{equation*}
\begin{split}
u'&=f_1 f_2\cdots f_{\ell-1}(u) = 
 [1..\ell \mid \ol{1}] \otimes [1..a \mid \ol{a}..\ol{b}]\\
u''&=f_2 f_3 \cdots f_\ell(u')=
 [1 \mid 3..\ell+1 \mid \ol{1}]\otimes [1..a \mid \ol{a}..\ol{b}]\\
u'''&=f_0(u'')=[1..\ell+1]\otimes [1..a \mid\ol{a}..\ol{b}].
\end{split}
\end{equation*}
By induction hypothesis, $v'''=R(u''')=[1..a \mid \ol{a}..\ol{b}]\otimes [1..\ell+1]$
and
\begin{equation*}
\begin{split}
v''&=e_0(v''')=[1..a-1\mid \ol{a-1}..\ol{b}\mid \ol{2}\ol{1}]\otimes [1..\ell+1]\\
v'&=e_\ell e_{\ell-1}\cdots e_2(v'')=[1..a-1\mid \ol{a-1}..\ol{b}\mid \ol{\ell+1}\mid\ol{1}]
 \otimes [1..\ell+1]\\
v&=e_{\ell-1} e_{\ell-2}\cdots e_1(v')=[1..a-1\mid \ol{a-1}..\ol{b}\mid \ol{\ell+1}\;\ol{\ell}]
 \otimes [1..\ell+1]
\end{split}
\end{equation*}
which by \eqref{eq:R crystal iso} and \eqref{eq:e-f} indeed yields 
$v=R(u)=[1..a-d\mid \ol{a-d}..\ol{b}\mid \ol{k'}..\ol{r}]\otimes [1..k']$.
Note also that $H(u''')=0$ by induction hypothesis. Since $f_0$ acts on
the left tensor factor in both $u''$ and $v''$, we have $H(u)=1$.

Next consider $d>1$ so that $u=[1..\ell\mid \ol{\ell}..\ol{r}]\otimes [1..a\mid \ol{a}..\ol{b}]$
with $r<\ell$. Then $r\le \ell-1 \le k'-2<b-2$ so that
\begin{equation*}
\begin{split}
u'&=f_2 f_3 \cdots f_r f_1 f_2\cdots f_{r-1}(u)=
 [1..\ell \mid \ol{\ell}..\ol{r+2} \mid \ol{2}\;\ol{1}] \otimes [1..a \mid \ol{a}..\ol{b}]\\
u''&=f_0(u')=
 [1..\ell+1\mid \ol{\ell+1}..\ol{r+2}]\otimes [1..a \mid \ol{a}..\ol{b}].
\end{split}
\end{equation*}
By induction hypothesis, $v''=R(u'')=[1..a-d+1 \mid \ol{a-d+1}..\ol{b}\mid \ol{k'}..\ol{r+2}]
\otimes [1..k']$ and 
\begin{equation*}
\begin{split}
v'&=e_0(v'')=[1..a-d \mid \ol{a-d}..\ol{b}\mid \ol{k'}..\ol{r+2}\mid \ol{2}\;\ol{1}]
 \otimes [1..k']\\
v&=e_{r-1} e_{r-2}\cdots e_1 e_r e_{r-1}\cdots e_2(v')=
 [1..a-d \mid \ol{a-d}..\ol{b} \mid \ol{k'}..\ol{r}]\otimes [1..k']
\end{split}
\end{equation*}
which yields $v=R(u)=[1..a-d \mid \ol{a-d}..\ol{b} \mid \ol{k'}..\ol{r}]\otimes [1..k']$ 
by \eqref{eq:R crystal iso} and \eqref{eq:e-f}. By induction hypothesis
$H(u'')=k'-\ell-1$. Since $f_0$ acts on the left tensor factor in both $u'$ and $v'$,
we have $H(u)=k'-\ell$ as claimed.

\newpage

\rotatebox{90}{$
\begin{array}{|c|l|l|l|} \hline
\text{Case} & R(u'\otimes u) & \text{conditions} & H(u'\otimes u)\\ 
\hline &&&\\[-2mm]
1 
& \left[1..\ell+1 \mid k'+1..k-1 \mid \ol{k-1}..\ol{q} \mid 
\ol{\ell+1}..\ol{r}\right] \otimes \left[1..k'\right]
& \text{$p=k$ and $k-q$ even}
& k'-\ell\\[2mm]
& \left[1..\ell \mid k'+1..p \mid \ol{p}..\ol{q} \mid \ol{\ell}..\ol{r} \right]
\otimes \left[1..k'\right]
& \text{otherwise}
& k'-\ell\\[2mm]
\hline &&&\\[-2mm]
2a
& [1..a-d\mid \ol{a-d}..\ol{b} \mid \ol{k'}..\ol{r}] \otimes
[1..k']
& k'<b\le a-d
& k'-\ell\\[2mm]
& [1..a-b+\ell+1 \mid k'+1..b-1 \mid \ol{a-b+\ell+1}..\ol{r}] \otimes
[1..k']
& k'<b; a-d<b
& k'-\ell\\[2mm]
& [1..\frac{k-k'+2\ell}{2} \mid \ol{\frac{k-k'+2\ell}{2}}..\ol{r}] \otimes
[1..\frac{k'+b-1}{2} \mid \ol{\frac{k'+b-1}{2}}..\ol{b}]
& \text{$b\le k'$, $k+k'$ even}
& k'-\frac{b+r-2}{2}\\[2mm]
& [1..b-1 \mid b+1..\frac{k-k'+2\ell+1}{2} \mid 
\ol{\frac{k-k'+2\ell+1}{2}}..\ol{r}]\otimes 
[1..\frac{k'+b}{2} \mid \ol{\frac{k'+b}{2}}..\ol{b+1}]
& \text{$b\le k'$, $k+k'$ odd}
& k'-\frac{b+r-3}{2}\\[2mm]
\hline &&&\\[-2mm]
2b
& [1..\frac{k-k'+2\ell}{2} \mid \ol{\frac{k-k'+2\ell}{2}}..\ol{r} \mid \ol{b-1}..\ol{q}]
\otimes [1..\frac{k'+b-1}{2} \mid \ol{\frac{k'+b-1}{2}}..\ol{b}]
& \text{$k+k'$ even}
& k'-\frac{r+q-2}{2}\\[2mm]
& [1..b-1 \mid b+1..\frac{k-k'+2\ell+1}{2} \mid \ol{\frac{k-k'+2\ell+1}{2}}..\ol{r} \mid 
\ol{b-1}..\ol{q}] \otimes [1..\frac{k'+b}{2} \mid \ol{\frac{k'+b}{2}}..\ol{b}]
& \text{$k+k'$ odd}
& k'-\frac{r+q-3}{2}\\[2mm]
\hline &&&\\[-2mm]
3
& [1..\ell \mid k'-k+b..p \mid \ol{\ell}..\ol{r}]
\otimes [1..k'-k+a \mid \ol{k'-k+a}..\ol{k'-k+b}]
& \text{$k<p$}
& k'-k+a-\ell\\[2mm]
& [1..\frac{k+p-2d}{2} \mid \ol{\frac{k+p-2d}{2}}..\ol{p+1} \mid \ol{\ell+d}..\ol{r}]
\otimes [1..\frac{k'+\ell+d}{2} \mid \ol{\frac{k'+\ell+d}{2}}..\ol{\ell+d+1}]
& \text{$k'\le p\le k$, $k+p$ even}
& \frac{k'-r+1}{2}\\[2mm]
& [1..\ell+d \mid \ell+d+2..\frac{k+p-2d+1}{2} \mid 
\ol{\frac{k+p-2d+1}{2}}..\ol{p+1} \mid \ol{\ell+d}..\ol{r}]
& \text{$k'\le p\le k$, $k+p$ odd}
& \frac{k'-r+2}{2}\\
& \otimes [1..\frac{k'+\ell+d+1}{2} \mid \ol{\frac{k'+\ell+d+1}{2}}..\ol{\ell+d+2}]
&&\\[2mm]
& [1..\frac{k+k'-2d}{2} \mid \ol{\frac{k+k'-2d}{2}}..\ol{p+1} \mid \ol{b-1}..\ol{r}]
\otimes [1..\frac{k'+b-1}{2} \mid \ol{\frac{k'+b-1}{2}}..\ol{b}]
& \text{$p<k'$, $k+k'$ even}
& \frac{k'-b+2d+1}{2}\\[2mm]
& [1..b-1 \mid b+1..\frac{k+k'-2d+1}{2} \mid \ol{\frac{k+k'-2d+1}{2}}..\ol{p+1} \mid
\ol{b-1}..\ol{r}] \otimes [1..\frac{k'+b}{2} \mid \ol{\frac{k'+b}{2}}..\ol{b+1}]
& \text{$p<k'$, $k+k'$ odd}
& \frac{k'-b+2d+2}{2}\\[2mm]
\hline
\end{array}
$}

\newpage

\rotatebox{90}{$
\begin{array}{|c|l|l|l|} \hline
\text{Case} & R(u'\otimes u) & \text{conditions} & H(u'\otimes u)\\ 
\hline &&&\\[-2mm]
4
& [1..\ell \mid k'-k+b..p \mid \ol{p}..\ol{q} \mid \ol{\ell}..\ol{r}]
\otimes [1..k'-k+a \mid \ol{k'-k+a}..\ol{k'-k+b}]
& \text{$k<p$}
& k'-k+a-\ell\\[2mm]
& [1..\frac{2\ell+k-p+2}{2} \mid k'+b-p..p-1 \mid \ol{p-1}..\ol{q} \mid
\ol{\frac{2\ell+k-p+2}{2}}..\ol{r}]
& \text{$k-2d<p\le k$, $k'<p$}
& k'-\ell+\frac{b-p-1}{2}\\
& \otimes [1..\frac{2k'+b-p-1}{2} \mid \ol{\frac{2k'+b-p-1}{2}}..\ol{k'+b-p}]
& \text{$k+p$ even, $p-q$ even}
& \\[2mm]
& [1..\frac{2\ell+k-p+1}{2} \mid k'+b-p+1..p \mid \ol{p}..\ol{q}
\mid \ol{\frac{2\ell+k-p+1}{2}}..\ol{r}]
& \text{$k-2d<p\le k$, $k'<p$}
& k'-\ell+\frac{b-p}{2}\\
& \otimes [1..\frac{2k'+b-p}{2} \mid \ol{\frac{2k'+b-p}{2}}..\ol{k'+b-p+1}]
& \text{$k+p$ odd, $p-q$ even}
& \\[2mm]
& [1..\frac{2\ell+k-p}{2} \mid k'+b-p..p \mid \ol{p}..\ol{q} \mid
\ol{\frac{2\ell+k-p}{2}}..\ol{r}]
& \text{$k-2d<p\le k$, $k'<p$}
& k'-\ell+\frac{b-p-1}{2}\\
& \otimes [1..\frac{2k'+b-p-1}{2} \mid \ol{\frac{2k'+b-p-1}{2}}..\ol{k'+b-p}]
& \text{$k+p$ even, $p-q$ odd}
& \\[2mm]
& [1..\frac{2\ell+k-p-1}{2} \mid k'+b-p+1..p+1 \mid \ol{p+1}..\ol{q}
\mid \ol{\frac{2\ell+k-p-1}{2}}..\ol{r}]
& \text{$k-2d<p\le k$, $k'<p$}
& k'-\ell+\frac{b-p}{2}\\
& \otimes [1..\frac{2k'+b-p}{2} \mid \ol{\frac{2k'+b-p}{2}}..\ol{k'+b-p+1}]
& \text{$k+p$ odd, $p-q$ odd}
& \\[2mm]
& [1..\ell+1+d \mid k'+b-p..\frac{k+p-2d-2}{2} \mid \ol{\frac{k+p-2d-2}{2}}..\ol{q}
\mid \ol{\ell+1+d}..\ol{r}]
& \text{$k'<p\le k-2d$}
& k'-\ell+\frac{b-p-1}{2}\\
& \otimes [1..\frac{2k'+b-p-1}{2} \mid \ol{\frac{2k'+b-p-1}{2}}..\ol{k'+b-p}]
& \text{$k+p$ even, $p-q$ even}
& \\[2mm]
& [1..\ell+1+d \mid k'+b-p+1..\frac{k+p-2d-1}{2} \mid \ol{\frac{k+p-2d-1}{2}}..\ol{q}
\mid \ol{\ell+1+d}..\ol{r}]
& \text{$k'<p\le k-2d$}
& k'-\ell+\frac{b-p}{2}\\
& \otimes [1..\frac{2k'+b-p}{2} \mid \ol{\frac{2k'+b-p}{2}}..\ol{k'+b-p+1}]
& \text{$k+p$ odd, $p-q$ even}
& \\[2mm]
& [1..\ell+d \mid k'+b-p..\frac{k+p-2d}{2} \mid \ol{\frac{k+p-2d}{2}}..\ol{q}
\mid \ol{\ell+d}..\ol{r}]
& \text{$k'<p\le k-2d$}
& k'-\ell+\frac{b-p-1}{2}\\
& \otimes [1..\frac{2k'+b-p-1}{2} \mid \ol{\frac{2k'+b-p-1}{2}}..\ol{k'+b-p}]
& \text{$k+p$ even, $p-q$ odd}
& \\[2mm]
& [1..\ell+d \mid k'+b-p+1..\frac{k+p-2d+1}{2} \mid \ol{\frac{k+p-2d+1}{2}}..\ol{q}
\mid \ol{\ell+d}..\ol{r}]
& \text{$k'<p\le k-2d$}
& k'-\ell+\frac{b-p}{2}\\
& \otimes [1..\frac{2k'+b-p}{2} \mid \ol{\frac{2k'+b-p}{2}}..\ol{k'+b-p+1}]
& \text{$k+p$ odd, $p-q$ odd}
& \\[2mm]
\hline
\end{array}
$}

\newpage

\rotatebox{90}{$
\begin{array}{|c|l|l|l|} \hline
\text{Case} & R(u'\otimes u) & \text{conditions} & H(u'\otimes u)\\ 
\hline &&&\\[-2mm]
4
& [1..\frac{2\ell+k-k'}{2} \mid b..p \mid \ol{p}..\ol{q}
\mid \ol{\frac{2\ell+k-k'}{2}}..\ol{r}]
\otimes [1..\frac{k'+b-1}{2} \mid \ol{\frac{k'+b-1}{2}}..\ol{b}]
& \text{$\frac{k+k'-2d}{2}< p\le k'$, $k+k'$ even}
& \begin{cases} \frac{2k'-r-q+4}{2} & \text{$p-q$ even}\\
                \frac{2k'-r-q+2}{2} & \text{$p-q$ odd}
\end{cases}\\[2mm]
\text{cont.}
&[1..\frac{2\ell+k-k'+1}{2} \mid b-1..p-1 \mid \ol{p-1}..\ol{q}
\mid \ol{\frac{2\ell+k-k'+1}{2}}..\ol{r}]
& \text{$\frac{k+k'-2d+1}{2}<p\le k'$}
& \frac{2k'-r-q+3}{2}\\
& \otimes [1..\frac{k'+b-2}{2} \mid \ol{\frac{k'+b-2}{2}}..\ol{b-1}]
& \text{$k+k'$ odd, $p-q$ even}
& \\[2mm]
& [1..\frac{2\ell+k-k'-1}{2} \mid b+1..p+1 \mid \ol{p+1}..\ol{q}
\mid \ol{\frac{2\ell+k-k'-1}{2}}..\ol{r}]
& \text{$\frac{k+k'-2d-1}{2}<p\le k'$}
& \frac{2k'-r-q+3}{2}\\
& \otimes [1..\frac{k'+b}{2} \mid \ol{\frac{k'+b}{2}}..\ol{b+1}]
& \text{$k+k'$ odd, $p-q$ odd}
& \\[2mm]
& [1..\ell+p-k'+d \mid b..\frac{k+k'-2d}{2} \mid \ol{\frac{k+k'-2d}{2}}..\ol{q}
\mid \ol{\ell+p-k'+d}..\ol{r}]
& \text{$k'<p+d\le \frac{k+k'}{2}$, $p\le k'$}
& \begin{cases} \frac{2k'-r-q+4}{2} & \text{$p-q$ even}\\
                \frac{2k'-r-q+2}{2} & \text{$p-q$ odd}
\end{cases}\\
& \otimes [1..\frac{k'+b-1}{2} \mid \ol{\frac{k'+b-1}{2}}..\ol{b}]
& \text{$k+k'$ even}
& \\[2mm]
& [1..\ell+p-k'+d \mid b-1..\frac{k+k'-2d-1}{2} \mid \ol{\frac{k+k'-2d-1}{2}}..\ol{q}
\mid \ol{\ell+p-k'+d}..\ol{r}]
& \text{$k'<p+d\le \frac{k+k'+1}{2}$, $p\le k'$}
& \frac{2k'-r-q+3}{2}\\
& \otimes [1..\frac{k'+b-2}{2} \mid \ol{\frac{k'+b-2}{2}}..\ol{b-1}]
& \text{$k+k'$ odd, $p-q$ even}
& \\[2mm]
& [1..\ell+p-k'+d \mid b+1..\frac{k+k'-2d+1}{2} \mid \ol{\frac{k+k'-2d+1}{2}}..\ol{q}
\mid \ol{\ell+p-k'+d}..\ol{r}]
& \text{$k'<p+d\le \frac{k+k'-1}{2}$, $p\le k'$}
& \frac{2k'-r-q+3}{2}\\
& \otimes [1..\frac{k'+b}{2} \mid \ol{\frac{k'+b}{2}}..\ol{b+1}]
& \text{$k+k'$ odd, $p-q$ odd}
& \\[2mm]
& [1..\ell \mid b..\frac{k+2p-k'}{2} \mid \ol{\frac{k+2p-k'}{2}}..\ol{q}
\mid \ol{\ell}..\ol{r}]
\otimes [1..\frac{k'+b-1}{2} \mid \ol{\frac{k'+b-1}{2}}..\ol{b}]
& \text{$p\le k'-d$, $k+k'$ even}
& \begin{cases} \frac{2k'-r-q+4}{2} & \text{$b-\ell$ even}\\
                \frac{2k'-r-q+2}{2} & \text{$b-\ell$ odd}
\end{cases}\\[2mm]
& [1..\ell \mid b-1..\frac{k+2p-k'-1}{2} \mid \ol{\frac{k+2p-k'-1}{2}}..\ol{q}
\mid \ol{\ell}..\ol{r}]
& \text{$p\le k'-d$}
& \frac{2k'-r-q+3}{2}\\
& \otimes [1..\frac{k'+b-2}{2} \mid \ol{\frac{k'+b-2}{2}}..\ol{b-1}]
& \text{$k+k'$ odd, $b-\ell$ even}
& \\[2mm]
& [1..\ell \mid b+1..\frac{k+2p-k'+1}{2} \mid \ol{\frac{k+2p-k'+1}{2}}..\ol{q}
\mid \ol{\ell}..\ol{r}]
& \text{$p\le k'-d$}
& \frac{2k'-r-q+3}{2}\\
& \otimes [1..\frac{k'+b}{2} \mid \ol{\frac{k'+b}{2}}..\ol{b+1}]
& \text{$k+k'$ odd, $b-\ell$ odd}
& \\[2mm]
\hline
\end{array}
$}

\newpage 

\rotatebox{90}{$
\begin{array}{|c|l|l|l|} \hline
\text{Case} & R(u'\otimes u) & \text{conditions} & H(u'\otimes u)\\ 
\hline &&&\\[-2mm]
5a
& [1..\ell \mid k'+1..a \mid \ol{a}..\ol{q}] \otimes [1..k']
& \text{$q>\ell+1$, $k'-\ell$ even}
& k'-\ell\\[2mm]
& [1..\ell+1 \mid k'+1..a-1 \mid \ol{a-1}..\ol{q} \mid \ol{\ell+1}] \otimes
[1..k']
& \text{$q>\ell+1$, $k'-\ell$ odd}
& k'-\ell\\[2mm]
& [1..\ell \mid b..\frac{k+2b-k'-2}{2} \mid \ol{\frac{k+2b-k'-2}{2}}..\ol{q}] 
\otimes [1..\frac{k'+b-1}{2} \mid \ol{\frac{k'+b-1}{2}}..\ol{b}]
& \text{$q\le \ell+1<b$, $k+k'$ even}
& \begin{cases} \frac{2k'-\ell-q+3}{2} & \text{$b-\ell$ even}\\
               \frac{2k'-\ell-q+1}{2} & \text{$b-\ell$ odd}
\end{cases}\\[2mm]
& [1..\ell \mid b-1..\frac{k+2b-k'-3}{2} \mid \ol{\frac{k+2b-k'-3}{2}}..\ol{q}] 
\otimes [1..\frac{k'+b-2}{2} \mid \ol{\frac{k'+b-2}{2}}..\ol{b-1}]
& \text{$q\le \ell+1<b$, $k+k'$ odd}
& \frac{2k'-\ell-q+2}{2}\\
&& \text{$b-\ell$ even}&\\[2mm]
& [1..\ell \mid b+1..\frac{k+2b-k'-1}{2} \mid \ol{\frac{k+2b-k'-1}{2}}..\ol{q}] 
\otimes [1..\frac{k'+b}{2} \mid \ol{\frac{k'+b}{2}}..\ol{b+1}]
& \text{$q\le \ell+1<b$, $k+k'$ odd}
& \frac{2k'-\ell-q+2}{2}\\
&& \text{$b-\ell$ odd}&\\[2mm]
& [1..\frac{k+2\ell-k'}{2} \mid \ol{\frac{k+2\ell-k'}{2}}..\ol{\ell+1}
\mid \ol{b-1}..\ol{q}]
\otimes [1..\frac{k'+b-1}{2} \mid \ol{\frac{k'+b-1}{2}}..\ol{b}]
& \text{$q,b\le \ell+1$, $k+k'$ even}
& \frac{k'+b-2q+1}{2}\\[2mm]
& [1..b-1 \mid b+1..\frac{k+2\ell-k'+1}{2} \mid 
\ol{\frac{k+2\ell-k'+1}{2}}..\ol{\ell+1} \mid \ol{b-1}..\ol{q}]
& \text{$q,b\le \ell+1$, $k+k'$ odd}
& \frac{k'+b-2q+2}{2}\\
& \otimes [1..\frac{k'+b}{2} \mid \ol{\frac{k'+b}{2}}..\ol{b+1}]
&& \\[2mm]
\hline
\end{array}
$}

\newpage 

\rotatebox{90}{$
\begin{array}{|c|l|l|l|} \hline
\text{Case} & R(u'\otimes u) & \text{conditions} & H(u'\otimes u)\\ 
\hline &&&\\[-2mm]
5b
& [1..\ell+d+1 \mid k'+1..a-1-d \mid \ol{a-1-d}..\ol{q} \mid 
 \ol{\ell+d+1}..\ol{r}]
& \text{$q\ge \ell+d+2$, $k'-r$ even}
& k'-\ell\\
& \otimes [1..k']
& \text{$a-b\ge d-1$}
& \\[2mm]
& [1..\ell+d \mid k'+1..a-d \mid \ol{a-d}..\ol{q} \mid 
 \ol{\ell+d}..\ol{r}]
& \text{$q\ge \ell+d+2$, $k'-r$ odd}
& k'-\ell\\
& \otimes [1..k']
& \text{$a-b\ge d-1$}
& \\[2mm]
& [1..\ell+a-b+2 \mid k'+1..b-2 \mid \ol{b-2}..\ol{q} \mid 
 \ol{\ell+a-b+2}..\ol{r}]
& \text{$q\ge \ell+d+2$, $k'-r$ even}
& k'-\ell\\
& \otimes [1..k'] 
& a-b<d-1
& \\[2mm]
& [1..\ell+a-b+1 \mid k'+1..b-1 \mid \ol{b-1}..\ol{q} \mid 
 \ol{\ell+a-b+1}..\ol{r}]
& \text{$q\ge \ell+d+2$, $k'-r$ odd}
& k'-\ell\\
& \otimes [1..k'] 
& a-b<d-1
& \\[2mm]
& [1..q-1 \mid b..a-\frac{q-r}{2} \mid \ol{a-\frac{q-r}{2}}..\ol{r}]
\otimes [1..\frac{k'+b-1}{2} \mid \ol{\frac{k'+b-1}{2}}..\ol{b}]
& \ell+d+1\ge q>\ell+1
& \begin{cases} \frac{2k'-r-q+4}{2} & \text{$k'-r$ even}\\
                \frac{2k'-r-q+2}{2} & \text{$k'-r$ odd}
\end{cases}\\
&& \text{$r-q$ even, $a-b\ge \frac{q-r-2}{2}$}&\\[2mm]
& [1..a-b+\frac{q+r}{2} \mid \ol{b-1}..\ol{q} \mid 
\ol{a-b+\frac{q+r}{2}}..\ol{r}]
\otimes [1..\frac{k'+b-1}{2} \mid \ol{\frac{k'+b-1}{2}}..\ol{b}]
& \ell+d+1\ge q>\ell+1
& \begin{cases} \frac{2k'-r-q+4}{2} & \text{$k'-r$ even}\\
                \frac{2k'-r-q+2}{2} & \text{$k'-r$ odd}
\end{cases}\\
&& \text{$r-q$ even, $a-b<\frac{q-r-2}{2}$}&\\[2mm]
&[1..q-1 \mid b-1..a-\frac{q-r+1}{2} \mid \ol{a-\frac{q-r+1}{2}}..\ol{r}]
& \text{$\ell+d+1\ge q>\ell+1$, $r-q$ odd}
& \frac{2k'-r-q+3}{2}\\
& \otimes [1..\frac{k'+b-2}{2} \mid \ol{\frac{k'+b-2}{2}}..\ol{b-1}]
& \text{$k'-r$ even, $a-b\ge \frac{q-r-3}{2}$}
& \\[2mm]
& [1..a-b+\frac{q+r+1}{2} \mid \ol{b-2}..\ol{q} \mid 
\ol{a-b+\frac{q+r+1}{2}}..\ol{r}]
& \text{$\ell+d+1\ge q>\ell+1$, $r-q$ odd}
& \frac{2k'-r-q+3}{2}\\
& \otimes [1..\frac{k'+b-2}{2} \mid \ol{\frac{k'+b-2}{2}}..\ol{b-1}]
& \text{$k'-r$ even, $a-b<\frac{q-r-3}{2}$}
& \\[2mm]
& [1..q-1 \mid b+1..a-\frac{q-r-1}{2} \mid \ol{a-\frac{q-r-1}{2}}..\ol{r}]
& \text{$\ell+d+1\ge q>\ell+1$, $r-q$ odd}
& \frac{2k'-r-q+3}{2}\\
& \otimes [1..\frac{k'+b}{2} \mid \ol{\frac{k'+b}{2}}..\ol{b+1}]
& \text{$k'-r$ odd, $a-b\ge \frac{q-r-1}{2}$}
& \\[2mm]
& [1..a-b+\frac{q+r-1}{2} \mid \ol{b}..\ol{q} \mid 
\ol{a-b+\frac{q+r-1}{2}}..\ol{r}]
& \text{$\ell+d+1\ge q>\ell+1$, $r-q$ odd}
& \frac{2k'-r-q+3}{2}\\
& \otimes [1..\frac{k'+b}{2} \mid \ol{\frac{k'+b}{2}}..\ol{b+1}]
& \text{$k'-r$ odd, $a-b<\frac{q-r-1}{2}$}
& \\[2mm]
\hline
\end{array}
$}

\end{document}